\documentclass{article}
\usepackage{geometry}
\usepackage{graphicx}
\usepackage{amsmath, amssymb, amsfonts}
\usepackage{algorithm, algpseudocode}
\usepackage{bm}
\usepackage{tikz}
\usepackage{stmaryrd}
\usepackage{appendix}
\usepackage{lineno}
\usepackage[affil-it]{authblk}
\usepackage{xcolor}
\usepackage{ulem}
\newcommand{\stkout}[1]{\ifmmode\text{\sout{\ensuremath{#1}}}\else\sout{#1}\fi}
\providecommand{\keywords}[1]
{
  \small	
  \textbf{\textit{Keywords---}} #1
}

\title{Exploring forms of the moist shallow water equations using a new compatible finite element discretisation}

\author[1,*]{Nell Hartney}
\author[2]{Thomas M. Bendall}
\author[1]{Jemma Shipton}

\affil[1]{Department of Mathematics and Statistics, University of Exeter, Exeter, UK}
\affil[2]{Dynamics Research, Met Office, Exeter, UK}
\affil[*]{Corresponding author: Nell Hartney, nh491@exeter.ac.uk}

\begin{document}

\maketitle

\begin{abstract}
    The moist shallow water equations offer a promising route for advancing understanding of the coupling of physical parametrisations and dynamics in numerical atmospheric models, an issue known as `physics-dynamics coupling'. Without moist physics, the traditional shallow water equations are a simplified form of the atmospheric equations of motion and so are computationally cheap, but retain many relevant dynamical features of the atmosphere. Introducing physics into the shallow water model in the form of moisture provides a tool to experiment with numerical techniques for physics-dynamics coupling in a simple dynamical model. In this paper, we compare some of the different moist shallow water models by writing them in a general formulation. The general formulation encompasses three existing forms of the moist shallow water equations and also a fourth, previously unexplored formulation. The equations are coupled to a three-state moist physics scheme that interacts with the resolved flow through source terms and produces two-way physics-dynamics feedback. We present a new compatible finite element discretisation of the equations and apply it to the different formulations of the moist shallow water equations in three test cases. The results show that the models capture generation of cloud and rain and physics-dynamics interactions, and demonstrate some differences between moist shallow water formulations and the implications of these different modelling choices.
\end{abstract}
\section{Introduction}

There are a number of significant processes in a weather or climate model that are not captured by the fluid dynamical equations that describe the behaviour of the atmosphere. These unresolved physical processes are parametrised, and often just called `physics', in contrast to the discretisation of the resolved flow which is known as the `dynamics'. The physics includes processes that occur on scales that are too small or fast to be resolved by the solution of the fluid equations, and non-fluid dynamical processes such as radiation. In a forecast model the physics must be incorporated with the fluid mechanics solution to allow these processes to influence the large-scale dynamics correctly. This is collectively known as physics-dynamics coupling, and, as described in, for example \cite{PDCreview}, the implementation of this coupling can have a significant impact on the model as well as presenting challenges for numerical methods. For instance, physics parametrisations can generate small-scale or grid-scale variability (and/or short-time variability), and it is important for numerical methods to be able to handle this appropriately. A common feature of physics schemes is the presence of discontinuities or switch-like processes which are also challenging from a numerical point of view.

We are interested in the potential that a simplified model has to investigate some of these numerical challenges associated with physics-dynamics coupling. Our approach is to use a model with simplified dynamics, but coupled to a physically-realistic moist physics scheme: the moist shallow water model, which couples moist physics interactions to the traditional shallow water dynamics.

The choice of using the shallow water equations is motivated by the fact that they are a simplification to the full three-dimensional atmospheric equations. They are computationally cheap yet still represent pertinent features of the atmosphere, such as Rossby waves and geostrophic adjustment. They are the simplest systems in geophysical fluid dynamics that exhibit separation between slow balanced solutions and fast unbalanced inertia-gravity waves, which make them a highly useful tool for numerical development. The addition to the shallow water equations of moisture as a sub-grid-scale process introduces complexities such as new timescales for moist physics processes and non-linear phase transitions which challenge numerical schemes. The fact that these complexities are not unique to moist physics but rather are relevant to other physics parametrisations too makes the moist shallow water equations a particularly good tool for investigating numerics questions about how the physics and dynamics processes interact with one another in the model. 

The idea of including moist physics processes in the shallow water equations is one that has been garnering some interest, with the primary motivation prompting the development of moist shallow water models being their use in modelling meteorological phenomena. This is because the highly significant role of water vapour in the Earth's climate makes models that include moist physics attractive from a meteorological point of view. For example, energy from phase changes can drive large-scale atmospheric motions and climatic effects, such as El Ni\~no, the Madden-Julian oscillation, the Hadley circulation, convectively-coupled gravity waves and monsoons \cite{StechmannMajda2006}. Investigations into these features have regularly turned to moist shallow water models as a framework to model some of these effects, offering as they do simple and inexpensive dynamics coupled to moist physics. Different frameworks for the inclusion of moisture in the shallow water system have been proposed by different authors (see, for example \cite{Gill1982, Bouchut2009, Lambaertsetal2011baroclinic, KLZ2020b, ZA2015, Yang2021, WC2014}), and these models have been used in studies of various meteorological effects (see, for example, \cite{Lambaertsetal2011jet, LahayeZeitlin2016, RostamiZeitlin2017, RostamiZeitlin2018improved, RostamiZeitlin2019, RostamiZeitlin2020, RostamiZeitlin2020MJO, RostamiZeitlin2021, RostamiZeitlin2022, RostamiZeitlin2022jets, VP2020}).

In spite of this extensive use of moist shallow water models, there is no comparison of the different models or investigation of the implications of various modelling choices. An understanding of the different modelling choices and a suite of test cases that can be run with different moist shallow water models is a necessary step to be able to use these models to investigate questions around physics-dynamics coupling.

This paper aims to address this gap by advancing model comparison in moist shallow water models, so as to aid progress in their use for physics-dynamics coupling investigations. We write the moist shallow water equations in a flexible form from which four different moist shallow water frameworks - including one previously unexplored - can be recovered, and present a discretisation of these equations using compatible finite elements, which is a first in moist shallow water modelling. The finite element method has become popular in recent times as a spatial discretisation approach in dynamical cores for a number of reasons, as described by Cotter in \cite{Cotter2023}. Finite element methods do not rely on a structured underlying mesh and are suitable on non-orthogonal grids, such as cubed-sphere and icosahedral grids. These quasi-uniform meshes of the sphere provide a solution to the issue of parallel scalability bottlenecks due to the convergence of grid lines around the poles in a latitude-longitude grid, described in, for example \cite{Lawrence2018}. Compatible finite element methods allow the preservation of discrete vector calculus identities (analogous to the C-grid for finite difference schemes) at the level of the finite element spaces and the mappings between them. This is key to ensuring the desirable numerical properties outlined by Staniforth and Thuburn in \cite{StaniforthThuburn2012}, such as supporting steady geostrophic modes, avoiding spurious mode branches, and offering the potential for schemes that will preserve energy and mass.

This has prompted the use of the compatible finite element method for the basis of the Met Office's next-generation dynamical core as described by Melvin \textit{et al.} in \cite{GungHoCartesian}. Our use of the compatible finite element method here is motivated by a desire to learn about physics-dynamics coupling in a way that will be relevant to this new dynamical core. We will outline our new discretisation of the flexible moist shallow water model using compatible finite elements and will demonstrate results in each of the moist shallow water formulations from three simple test cases (motivated by standard shallow water test cases in the literature): steady-state geostrophically balanced flow (based on test 2 from the Williamson \textit{et al.} test suite \cite{Williamson1992}), flow over an isolated mountain (based on test 5 from the same suite) and the unstable jet of Galewsky et al \cite{Galewsky2004}. Each of our tests will use a semi-implicit quasi-Newton time-stepping scheme in the style of that used by the Met Office's ENDGame \cite{Wood2014inherently} and GungHo \cite{GungHoCartesian} dynamical cores.

The novelty of this paper lies in the use of compatible finite elements to discretise this model, a technique that, to our knowledge, has not been used in moist shallow water modelling before. Our discretisation is achieved using a generalised form of the moist shallow water equations, a form which also facilitates a fourth, previously unexplored, moist shallow water formulation. We also detail the set-ups for three test cases suitable for use with all four different formulations of the moist shallow water equations, and use these tests to uncover differences between formulations and implications of different modelling choices. The initial conditions for the tests are carefully derived from balance conditions and are detailed thoroughly in each moist shallow water framework so as to be easily reproducible. The idea is that these tests can act as reference solutions for other practitioners and can provide information around the implication of modelling choices in different moist shallow water frameworks to help inform modelling decisions, all with the aim of advancing model comparison for moist shallow water models in order to optimise the model as a tool for physics-dynamics coupling investigations.

The structure of the article is as follows: in Section 2 we give an overview of the three most commonly-used moist shallow water models and their use in the literature. Section 3 describes the general formulation for the moist shallow water models and how to recover the different formulations from the general one. This section thus provides us with the framework with which we can compare formulations. Both the compatible finite element discretisation of the model and our time-stepping scheme are dealt with in Section 4, which provides details of the new discretisation and demonstrates how our approach to solving the equations is highly relevant in a contemporary numerical weather prediction paradigm. Finally, the test cases and results in Section 5 offer representative solutions to compare with and information on modelling choices.

\section{Moist shallow water models}

\subsection{The shallow water equations}
The shallow water equations describe flows where horizontal length scales are much larger than vertical length scales. They are derived by integrating the rotating Euler equations between two material surfaces, before assuming that the horizontal flow is independent of vertical coordinate and applying the hydrostatic balance condition \cite{Zeitlinbook}. Written here with rotation and without any dissipative effects, the shallow water equations are given by:

\begin{align}
    &\frac{\partial{\bm{u}}}{\partial{t}} + (\bm{u} \cdot \bm{\nabla})\bm{u} + f \bm{\hat{k}} \times \bm{u}= -g \bm{\nabla} (D + B), \\
    &\frac{\partial{D}}{\partial{t}} + \bm{\nabla} \cdot (\bm{u}D) = 0,
\end{align}
where $\bm{u}$ is horizontal velocity, $D$ is the thickness of the fluid layer, $g$ is acceleration due to gravity, $f$ is the Coriolis parameter, $\bm{\hat{k}}$ is the radial unit vector and $B$ is the height of the bottom surface.

In this traditional form the shallow water equations are a purely dynamical model, with no representation of sub-grid-scale physical processes. To be relevant for physics-dynamics coupling the model must include some physics processes. Moisture-related sub-grid-scale processes are a natural choice, having been included in the shallow water equations in many different ways in the past. Approaches to incorporating moisture into the shallow water system fall broadly into two groups: models where there is no explicit moisture variable but the effects of moisture are included in the shallow water dynamics, and models that add one or more moisture variables that feedback on the existing shallow water dynamics.

In the first group are frameworks that model the effect of moisture in the shallow water system, rather than modelling moisture directly. The moist shallow water model of Yang and Ingersoll (\cite{Yang2013}, \cite{Yang2014} and \cite{Yang2021}) falls into this category. In this approach convection is parametrised as a short-duration forcing term in the depth equation, triggered locally by changes in layer thickness. A related idea is used by W{\"u}rsch and Craig in their model for use in development and testing of data assimilation algorithms \cite{WC2014}. In this model the definition of geopotential is altered in regions of decreased geopotential such that its gradient forces fluid into the region, simulating the convective effects of a cloud on geopotential. This model is adapted by Kent, Bokhove and Tobias in \cite{Kent2017} for use in data assimilation  by excluding the diffusion terms in the original model and including rotation, thus leading to a so-called 1.5-dimensional model. This version of the model is used as a testbed for investigating data assimilation algorithms in \cite{Kent2020}, and is modified in \cite{Bokhove2022} by adding a second passive layer and by replacing the uniform layer density with uniform potential temperature.

Our focus is on the second group of models; that is, models that add one or more equations for moist variables to the usual shallow water set and then couple the additional equations to the pre-existing ones in various ways. We have chosen this group of models because this explicit addition of moisture mimics what happens in dynamical cores, making them more suitable for physics-dynamics coupling investigations.

\subsection{The moist convective shallow water equations}
The earliest model in this style is the seminal work by Gill \cite{Gill1982}. Here an equation for a new moisture variable is added to the shallow water equations and convection is modelled via a sink in the depth equation. The model is used in a dam break problem, a Rossby adjustment problem, and to model wave propagation through moist and dry regions. The Gill model has also been used in studies of tropical data assimilation in, for example, \cite{Vzagar2008} and \cite{Vzagar2012}.  Similar models to the Gill framework were used in studies of precipitation fronts by Frierson et al. \cite{Frierson2004}, Stechman and Majda \cite{StechmannMajda2006}, and Pauluis et al. \cite{Pauluis2008}. These use a linear rotating shallow water-type model with a vertically averaged moisture field, with convection parameterised by a non-linear precipitation term. This non-linear source term gives rise to precipitation front solutions. Fronts are studied in the limit of infinitely fast relaxation in \cite{Frierson2004} and shown to be realisable in finite relaxation time in \cite{StechmannMajda2006}. Interactions between a precipitation front and an incident disturbance are studied in \cite{Pauluis2008}.

Also inspired by the Gill approach, Bouchut et al. \cite{Bouchut2009} outlined a systematic method to include moisture, condensation and related latent heat release in the shallow water equations. The authors call this the moist convective shallow water model as it adds an advection equation for a moisture variable to the classical shallow water set, and couples the moisture equation to the depth equation through a sink term. This model is used in Lambaerts et al. \cite{Lambaertsetal2011jet} to study the effect of moisture on the development of the Bickley jet barotropic instability. In the baroclinic version of the moist convective model proposed by Lambaerts et al. in \cite{Lambaertsetal2011baroclinic}, the model is extended to two layers, where the mass loss or sink term in the depth equation is modelled as mass exchange between the layers. This two-layer moist convective model is used to study the evolution of the Bickley jet baroclinic instability in \cite{Lambaerts2012}. Evaporation from a wet surface is included in the barotropic moist convective model used to study instabilities of hurricane-like vortices in \cite{LahayeZeitlin2016}, where it appears as a source term proportional to the intensity of the wind in the moisture equation.

In all of these moist convective shallow water models the condensed water vapour disappears from the model, when in reality condensed water in the atmosphere remains in the form of clouds, and precipitation is only triggered when water droplets reach a critical size. The improved moist convective shallow water model of Rostami and Zeitlin \cite{RostamiZeitlin2018improved} aims to overcome some of these drawbacks. The model extends the two-layer model of Lambaerts et al. \cite{Lambaertsetal2011baroclinic} by including, as well as condensation, an inverse vaporisation phase transition to allow the conversion between precipitable water and vapour. The model also includes a source of bulk moisture in the lower of the two layers due to surface evaporation. Topography is included in this improved moist-convective model in \cite{RostamiZeitlin2022jets}, where the model is used to investigate the evolution of waves on the background of a easterly jet that crosses a land-sea boundary.

Variations of these one- and two-layer moist convective models are used in a number of studies of tropical dynamics. These include investigations of the instabilities of tropical cyclone-like vortices (\cite{LahayeZeitlin2016}, \cite{RostamiZeitlin2017}, \cite{RostamiZeitlin2020} and \cite{RostamiZeitlin2022}, for example) and studies of equatorial modons and the Madden-Julian oscillation (\cite{RostamiZeitlin2019}, \cite{Zhao2021}, \cite{RostamiZeitlin2021}, and \cite{VP2020}, \cite{RostamiZeitlin2020MJO}, for example.) 

The formulation of the moist convective model we use is that of Bouchut et al. \cite{Bouchut2009}. Condensation is triggered when the moisture variable exceeds some saturation threshold and acts to lower the layer thickness, through a sink term proportional to the condensation in the depth equation. The model is known as moist convective because the moist response triggers a convective one, where $-\beta P$ represents a mass loss due to convection out of the shallow water layer:
\begin{align}
    &\frac{\partial{\bm{u}}}{\partial{t}} + (\bm{u} \cdot \bm{\nabla})\bm{u} + f \bm{\hat{k}} \times \bm{u}= -g \bm{\nabla} (D + B), \\
    &\frac{\partial{D}}{\partial{t}} + \bm{\nabla} \cdot (\bm{u}D) = -\beta P, \\
    & \frac{\partial{Q}}{\partial{t}} + \bm{\nabla} \cdot (\bm{u}Q) = -P,
\end{align}
where $Q$ represents the amount of water vapour in the layer (the specific humidity integrated over the layer), and so obeys a conservative form of the transport equation. Though topography $B$ is not included in the original model we include it here. Moisture satisfies the relaxation relation of the type used in general circulation models \cite{BettsMiller}:
    \begin{equation} \label{eq:BettsMiller}
        P = \frac{Q - Q_s}{\tau} H\left(Q - Q_s\right),
    \end{equation}
where $Q_s$ is a saturation value, $H(\cdot)$ is the Heaviside step function and $\tau$ is the timescale of condensation. The simplest parametrisation involves choosing a constant value $Q_s$, but more sophisticated parametrisations (such as a depth- or pressure-dependent saturation function) are also possible. In this model condensation and precipitation are synonymous, and condensed water vapour is not included in the dynamics.

As noted by Bouchut et al. in \cite{Bouchut2009}, the (limited) loss of mass through the upper layer of the model means that this moist convective shallow water model should be interpreted as the lower layer in a more complete multi-layer model. This idea is also addressed by Zeitlin \cite{Zeitlinbook}, where this single layer model is considered as part of a configuration of a two-layer model with a dry upper layer, in the limit where the upper layer is much deeper than the lower one.

\subsection{The moist convective thermal shallow water equations}
The thermal shallow water equations are a variant of the usual shallow water model that include the representation of thermodynamic processes. They are derived by vertically averaging the primitive equations as in the derivation of the usual shallow water system, but with the hypothesis of uniform temperature or density relaxed (see, for example \cite{Zeitlinbook}). The idea of incorporating moist effects into this thermal shallow water model was taken up by Kurganov, Liu and Zeitlin in \cite{KLZ2020} with their implementation of the moist convective thermal shallow water equations. The main advantage of this model over the moist convective one is that it allows the representation of both the mechanical and thermal effects of convection, by splitting the effect of the latent heat into two parts: heating, or local increase of potential temperature, and convective flux. This moist convective thermal shallow water model is used by the same authors in \cite{KLZ2021} to model tropical cyclone-like vortices over oceanic warm and cold pools, and includes an island-like topography. Rostami, Zhao and Petri extend the model to a two-layer one in \cite{Rostami2022}, where they present a possible equatorial adjustment scenario for the genesis and dynamics of the Madden-Julian oscillation.

We use the model in its single layer version in our formulation (written here with topography, $B$) where the variable $b$ represents a buoyancy field:
\begin{align}
    &\frac{\partial{\bm{u}}}{\partial{t}} + (\bm{u} \cdot \bm{\nabla})\bm{u} + f \bm{\hat{k}} \times \bm{u}= -b \bm{\nabla} (B + D) - \frac{D}{2} \bm{\nabla} b,  \\ 
    &\frac{\partial{D}}{\partial{t}} + \bm{\nabla} \cdot (\bm{u} D) = -\beta_1 P, \\
    &\frac{\partial b}{\partial t} + (\bm{u} \cdot \bm{\nabla})b = \beta_2 P, \\ 
    & \frac{\partial{Q}}{\partial{t}} + \nabla \cdot (\bm{u}Q) = -P.
\end{align}
Like the moist convective model, moisture satisfies a relaxation relation with a prescribed saturation function, and the single moisture variable means that condensation and precipitation are synonymous and that condensed water vapour is not included in the dynamics.

\subsection{The moist thermal shallow water equations}
The idea of modelling the heating effect of moisture in the shallow water system is also dealt with by Zerroukat and Allen in \cite{ZA2015}, though with a different approach. Rather than beginning with a version of the usual shallow water equations and adding moisture, they derive a moist shallow water system from the three-dimensional Boussinesq approximation of the hydrostatic Euler equations. Unlike the usual Boussinesq equations, density is allowed to vary slightly with temperature, meaning the addition to the standard Boussinesq equations of an advection equation for a temperature/density variable and the addition to the momentum equation of buoyancy-related terms. The moist shallow water equations are derived by integrating this moist Boussinesq set over the depth of a shallow fluid layer. Moisture exists in three forms; vapour, cloud water and rain, and a three-state phase change model provides two-way physics-dynamics feedback. The framework is used in this study to demonstrate various test problems in the style of the Williamson et al. test cases \cite{Williamson1992} (steady-state geostrophically balanced flow, flow over an isolated mountain, mid-latitude jets over two isolated mountains, and zonal flow over Earth orography). The model is also used by Santos and Peixoto in \cite{SantosPeixoto2021} in their study of the impact of local grid refinement on cloud and rain formation over the Andes mountains, and by Ferguson and Jablonowski in \cite{FergusonJablonowski2019} to test the ability of adaptive mesh refinement to track and resolve sharp gradients and small-scale flow filaments. Written here with topography $B$, the model is given by:

\begin{align}
    &\frac{\partial{\bm{u}}}{\partial{t}} + (\bm{u} \cdot \bm{\nabla})\bm{u} + f \bm{\hat{k}} \times \bm{u}= - g \bm{\nabla} (D + B) + g \theta \bm{\nabla} B + \frac{g}{D} \bm{\nabla} \left(\frac{1}{2} D^2 \theta \right), \label{eq:mt_sw_eqns_start} \\
    &\frac{\partial{D}}{\partial{t}} + \bm{\nabla} \cdot (\bm{u}D) = 0, \\
    &\frac{\partial{\theta}}{\partial{t}} + (\bm{u} \cdot \bm{\nabla})\theta = S_{\theta}, \\
    & \frac{\partial{q^{(k)}}}{\partial{t}} + (\bm{u} \cdot \bm{\nabla)}{q^{(k)}} = S_q^{(k)}. \label{eq:mt_sw_eqns_end}
\end{align}
The new variable $\theta$ represents the depth-averaged deviation from a background temperature. Here $q^{(k)}$ where $k \in 1, 2,3$ are the moisture fields (water vapour, cloud, and rain), and $S_q^{(k)}$ are source terms due to phase changes between the fields. In contrast to the variable $Q$ these moist fields are specific humidities, and so satisfy the advective form of the transport equation. $S_{\theta}$ is a source term for temperature due to latent heat release from moisture phase changes. Phase changes are described using a three-state moist physics scheme. Whenever the water vapour exceeds a threshold dictated by a prescribed saturation function (here a function of $D$ and $\theta$), a fraction of the excess is converted to cloud, with a corresponding latent heat release that heats the atmosphere locally. When the atmosphere is sub-saturated and there is cloud present a fraction of the cloud evaporates with a local cooling effect. When cloud passes a certain threshold a fraction of the excess is converted to rain, again with a corresponding latent heat release.  

A comparison between the moist convective equations, the moist convective thermal equations and the moist thermal equations is shown in Table \ref{table:literature_table}. 

\begin{table}
\begin{center}
\begin{tabular}{||c c c c c||} 
 \hline
 Framework & prognostic b & $S_D$ & $S_b$ & physics \\ [0.5ex] 
 \hline\hline
 moist convective & \text{\sffamily X}  & \checkmark & \text{\sffamily X} & one-way   \\
\hline
moist convective thermal & \checkmark & \checkmark & \checkmark & one-way \\
 \hline
moist thermal  & \checkmark & \text{\sffamily X} & \checkmark & 3-state
\\
\hline
\end{tabular}
\caption{A summary of the variations between the three most commonly-used moist shallow water models in the literature. The moist convective thermal model and the moist thermal model have prognostic thermal or buoyancy variables (denoted by $b$), while this variable can be considered to be constant in the moist convective model. $S_D$ is the right-hand side or source term in the depth equation, and $S_b$ the same for the buoyancy or potential temperature equation. These source terms are either present (\checkmark) or absent (\text{\sffamily X}) in each equation set. Both the moist convective model and the moist convective thermal model are usually coupled to a two-state physics scheme where conversions between water vapour and cloud are irreversible, called one-way physics in the table. The moist thermal model is coupled to a three-state moist physics scheme with water vapour, cloud and rain, where phases changes between water vapour and cloud are reversible.}
\label{table:literature_table}
\end{center}
\end{table}

\section{General moist shallow water formulation}

The moist convective model, the moist convective thermal model and the moist thermal model can be viewed as, respectively, a model where moist thermodynamics feeds back on (1) the depth equation, (2) the depth equation and the velocity equation, and (3) the velocity equation only. With this paradigm in mind, we consolidate these three approaches into one flexible model, with source terms that vary depending on the choice of moisture feedback. This requires writing both models that include thermal effects with the same thermal variable. We take $b = g(1 - \theta)$, which bears similarity to the relationship between the potential temperature and buoyancy in the Boussinesq equations in, for example, \cite{Zeitlinbook}, \cite{Vallis2019Book}, and \cite{Vallis2019RainyBenard} (noting that here $\theta$ is dimensionless, following Zerroukat and Allen \cite{ZA2015}), and we replace the temperature variable $\theta$ with $1 - b/g$ everywhere in the moist thermal shallow water equations (Equations \eqref{eq:mt_sw_eqns_start} - \eqref{eq:mt_sw_eqns_end}).

The advantage of defining this relationship between $\theta$ and $b$ is that the velocity equation in the moist thermal model will share a right-hand side with the moist convective thermal model. This means that all three models are encapsulated by one equation set, where different choices for the source terms and the buoyancy recover each of the individual models. Thus the equation set that we will discretise in the following section is given by: 

\begin{align}
    &\frac{\partial{\bm{u}}}{\partial{t}} + (\bm{u} \cdot \bm{\nabla})\bm{u} + f \bm{\hat{k}} \times \bm{u} = -b \bm{\nabla} (D + B) - \frac{D}{2} \bm{\nabla} b, \\
    &\frac{\partial{D}}{\partial{t}} + \bm{\nabla} \cdot (\bm{u}D) = \beta_1 S_D, \label{eq:general_D_eq} \\
    &\frac{\partial b}{\partial t} + (\bm{u} \cdot \bm{\nabla})b = \beta_2 S_b, \\
    & \frac{\partial{q^{(k)}}}{\partial{t}} + (\bm{u} \cdot \bm{\nabla}){q^{(k)}} = S_q^{(k)},
\end{align}
where $B$ is the topography, $S_D$ is a source term for the depth equation that captures the convective effect of phase changes and $S_b$ is a source term for the buoyancy equation that captures the heating effect of phase changes. The source terms $S_D$, $S_b$ and $S_q$ and the relationship between them come from the moist physics scheme. We have coupled the equations to two different moist physics schemes which are described below. The moisture variables $q^{(k)}$ represent mixing ratios, where $Q$ is the integral of $q$ over the depth of the shallow water layer, as in Zerroukat and Allen's moist thermal shallow water model \cite{ZA2015}. The mixing ratio variables obey the advective form of the transport equation, while the integrated variables Q would obey the conservative form of the transport equation. Moisture conservation could still be obtained using integrated variables through the inclusion of the depth field in the moisture transport equations, akin to the approach of \cite{FergusonJablonowski2019} in their moist shallow water formulation and the approach of \cite{Bendall2023trilemma} in the compressible Euler equations.

Writing the moist shallow water equations in a general form like this exposes a fourth moist shallow water formulation, not previously explored. With this idea, buoyancy is non-constant so that there is still a non-linear thermal contribution to the pressure gradient, but $\beta_2 = 0$ so that phase changes cause no latent heat impact and only a convective one. This formulation we call the \textit{moist convective pseudo-thermal} model, and is straightforward to implement in our discretisation since it simply involves switching on and off source terms in the depth and buoyancy equations. 

Table \ref{table:sources_table} outlines how each of the individual moist shallow water models (including the moist convective pseudo-thermal model) are recovered from the general model by making different choices for the buoyancy and the $\beta_1$ and $\beta_2$ parameters.

\begin{table}
\begin{center}
\begin{tabular}{||c c c c c||} 
 \hline
 Framework & $b$ & $\beta_1$ & $\beta_2$ & \\ [0.5ex] 
 \hline\hline
 moist convective & $g$  & non-zero & 0 & \\
\hline
moist convective thermal  & non-constant & non-zero & non-zero & \\
\hline
moist thermal & non-constant & 0 & non-zero &\\
\hline
moist convective pseudo-thermal  & non-constant & non-zero & 0 & \\
\hline
\end{tabular}
\caption{Choices for $b$, $\beta_1$ and $\beta_2$ in the flexible framework that recover each of the individual frameworks.}
\label{table:sources_table}
\end{center}
\end{table}

These general moist shallow water equations can be coupled to different physics schemes, and we have implemented both the one-way scheme described in Bouchut et al. \cite{Bouchut2009} and previously used with the moist convective and moist convective thermal shallow water equations, as well as the three-state physics scheme used with the moist thermal shallow water equations by Zerroukat and Allen \cite{ZA2015}.

The one-way physics scheme was coupled to the moist convective equations in \cite{Bouchut2009}, and to the moist convective thermal shallow water equations in \cite{KLZ2020b}. In this scheme moisture exists as water vapour which condenses when it exceeds some saturation value according to Equation \eqref{eq:BettsMiller}. To couple this scheme to our model we take: 
\begin{equation}
    S_v = -\frac{q_v - q_{sat}}{\tau} H(q_v-q_{sat}),
\end{equation}
where $S_v$ is the source term for vapour, $q_v$ represents the vapour and $q_{sat}$ is the prescribed saturation function. $S_D$ and $S_b$ are then equal to $S_v$, and $S_c = - S_v$. We define a saturation function which depends on depth and buoyancy in the models with a non-constant buoyancy field and depends on depth only in the moist convective model, which has a time-constant buoyancy field. In our implementation of this scheme we take $\tau$ to be the model timestep so that condensation is considered instantaneous. 

Our implementation of the three-state moist physics scheme follows that of Zerroukat and Allen \cite{ZA2015}. Moisture appears in three forms: vapour $q_v$, cloud $q_c$ and rain $q_r$. A saturation function $q_{sat}$ is prescribed in the same way as in the one-way physics scheme above (either depth-dependent or depth and buoyancy-dependent, depending on the moist shallow water framework). Whenever the local value of the vapour $q_v$ exceeds the saturation function then a portion of the excess $(q_v - q_{sat})$ will condense into cloud over a timescale $\tau_v$. This causes either a local heating of the atmosphere or a convective effect or both, depending on the moist shallow water framework. This conversion can be written as:
\begin{align}
    q_v > q_{sat} : C &= \gamma_v \frac{q_v - q_{sat}}{\tau_v}, \\
    &q_v \rightarrow q_v - \Delta t C, \label{eq:v_increment} \\ 
    &q_c \rightarrow q_c + \Delta t C, \\
    &D \rightarrow D - \Delta t \beta_1 C, \label{eq:D_increment} \\
    &b \rightarrow b - \Delta t \beta_2 C, \label{eq:b_increment}
\end{align}
where $\gamma_v$ is a conversion rate that depends on the depth and/or buoyancy and is discussed below.

In a similar way, in a sub-saturated atmosphere where cloud is present, cloud will evaporate:
\begin{align}
    q_v < q_{sat} : C &= \text{min}\left [\frac{q_c}{\Delta t}, \gamma_v \frac{q_{sat} - q_v}{\tau_v} \right ], \\
    &q_v \rightarrow q_v + \Delta t C, \\
    &q_c \rightarrow q_c - \Delta t C, \\
    &D \rightarrow D + \Delta t \beta_1 C, \\
    &b \rightarrow b + \Delta t \beta_2 C,
\end{align}
where the minimum function is used to ensure that only the available cloud can be evaporated.

If the saturation function varies strongly with buoyancy or depth, a phase change can cause a change in these fields that means that the atmosphere is immediately returned to its previous sub-saturated/over-saturated state, prompting another phase change which will return it again to the second state. This two-step oscillation is avoided by converting only a fraction of the cloud or vapour each time a change between these phases happens. This fraction, $\gamma_v$,  is derived below, and is chosen so that the final state is close to saturation at the new temperature.

Cloud is converted to rain when $q_c$ exceeds a threshold $q_{precip}$. A fraction of the excess is converted to rain over a timescale $\tau_r$ with a conversion rate $\gamma_r$, a constant since the process is not reversible and there is no danger of a two-step oscillation:
\begin{align}
    q_c > q_{precip} : P &=  \gamma_r \frac{q_c - q_{precip}}{\tau_r}, \\
    &q_r \rightarrow q_r + \Delta t P, \\
    &q_c \rightarrow q_c - \Delta t P,
\end{align}
Because there is no sink for $q_r$, it can be interpreted as a moisture phase produced from the phase change model. In our implementation of this physics scheme we take the timescales for both the conversion between water vapour and cloud $\tau_v$ and the conversion between cloud and rain $\tau_r$ to be the same as the timestep of the model, meaning that these processes are considered instantaneous. The rain does not feed back on the dynamics in either Zerroukat and Allen \cite{ZA2015} or in our implementation of the moist physics scheme, though in our implementation the rain is accumulated while in Zerroukat and Allen \cite{ZA2015} it is advected around.

Putting these three conversion processes together the source terms can be concisely written as:
\begin{equation}
\begin{split}
    \Delta q_v &= \text{max}\left [ 0, \gamma_v \left(\frac{q_v - q_{sat}}{\tau_v} \right)\right ], \\
    \Delta q_c &= \text{min}\left [ \frac{q_c}{\Delta t}, \text{max}\left[0, \gamma_v \left( \frac{q_{sat} - q_v}{\tau_v} \right) \right ]\right ], \\
    \Delta q_r &= \text{max}\left [ 0, \gamma_r \left(\frac{q_c - q_{precip}}{\tau_{r}} \right)\right ], \\
    S_v &= \Delta q_c - \Delta q_v, \\
    S_c &= \Delta q_v - \Delta q_c - \Delta q_r, \\
    S_r &= \Delta q_r, \\
    S_D &= \Delta q_c - \Delta q_v, \\
    S_b &= \Delta q_c - \Delta q_v .
\end{split}
\end{equation}

Coupling the three-state scheme to our model means that we have introduced additional moisture variables into both the moist convective and the moist convective thermal frameworks, which in their original formulations included only one moist species and considered condensation and precipitation synonymous. It also introduces to these frameworks the idea of a reversible conversion between water vapour and cloud, where cloud can evaporate in a sub-saturated atmosphere. This is an important aspect of physics-dynamics coupling in atmospheric models and its inclusion furthers our aim of a useful model for exploring physics-dynamics coupling questions. A comparison of this three-state physics scheme and the irreversible physics scheme is given in Appendix A, and the values of $\gamma_v$, $\beta_1$, $\beta_2$, $\gamma_r$ and $q_{precip}$ are given below.

Following Appendix A in Zerroukat and Allen \cite{ZA2015}, for the moist thermal frameworks (moist thermal, moist convective thermal and moist convective pseudo-thermal) the saturation function is dependent on both depth and buoyancy\footnote{We introduce a factor of $gH$ to the saturation function given in Zerroukat and Allen \cite{ZA2015} to give a dimensionless number, but note that we modify the scaling factor $q_0$ accordingly to give a saturation function the same size as that in \cite{ZA2015}.}:
\begin{equation}
    q_{sat} = \frac{q_0 H}{D + B} e^{20(1-b/g)}, \label{eq:sat_func_thermal}
\end{equation}
where $q_0$ is a scaling factor given below in the results section, $H$ is the background depth and $B$ is the topography.
For the moist convective tests the saturation function varies only with depth:
\begin{equation}
    q_{sat} = \frac{q_0 H }{D + B} e^{20 \theta}. \label{eq:sat_func_depth}
\end{equation}
The $\theta$ term here is a constant in time and can be thought of as a proxy for how temperature varies with latitude. For all moist convective tests we include this spatially-varying $\theta$ so as to keep the set-ups similar between formulations.

We derive $\gamma_v$ by considering the increment to vapour, following Equation \eqref{eq:v_increment} in the physics scheme:
\begin{equation}
    q_v^{n+1} = q^n_v - \Delta t C = q^n_v - \gamma_v\left[q_v^n - q_{sat}(D^{n},b^{n})\right],
\end{equation}
where $n$ indicates the time level and we take $\tau_v$ to be the timestep $\Delta t$. 
We consider that this final vapour should be at saturation at the next time level, such that we have
\begin{equation}
q_v^{n+1} = q_{sat}(D^{n+1},b^{n+1}),
\end{equation}
and
\begin{equation}
    q_{sat}(D^{n+1},b^{n+1}) = q^n_v - \gamma_v\left[q_v^n - q_{sat}(D^{n},b^{n})\right]. \label{eq:q_sat_sub_in}
\end{equation}
We also have, from Equations \eqref{eq:D_increment} and \eqref{eq:b_increment}, that:
\begin{subequations}
\begin{align}
    D^{n+1} = D^n - \Delta t\beta_1C &= D^n - \beta_1\gamma_v \left[q_v^n - q_{sat}(D^{n},b^{n})\right], \label{eq:D_n+1} \\ 
    b^{n+1} = b^n - \Delta t\beta_2C &= b^n - \beta_2\gamma_v \left[q_v^n - q_{sat}(D^{n},b^{n})\right]. \label{eq:b_n+1}
\end{align}
\end{subequations}
We make these substitutions in \eqref{eq:q_sat_sub_in} and then assume that:
\begin{equation}
|\beta_1\gamma_v \left[q_v^n - q_{sat}(D^{n},b^{n})\right]| \ll |D^n|
\quad \text{and} \quad
|\beta_2\gamma_v \left[q_v^n - q_{sat}(D^{n},b^{n})\right]| \ll |b^n|.
\end{equation}
We expand the left-hand side of Equation \eqref{eq:q_sat_sub_in} using a first-order Taylor expansion and rearrange to give:
\begin{equation}
    \gamma_v = \frac{1}{1 - \beta_1\frac{\partial q_{sat}}{ \partial D}(D^n,b^n) - \beta_2\frac{\partial q_{sat}}{ \partial b}(D^n,b^n)}.
\end{equation}
Returning to the saturation function \eqref{eq:sat_func_thermal}, we have
\begin{equation}
    \frac{\partial q_{sat}}{ \partial D}(D^n,b^n) = \frac{\partial}{\partial D} \left( \frac{q_0 H}{D + B} e^{20(1-b/g)}\right)
    = -\frac{1}{D + B} q_{sat}(D^n, b^n),
\end{equation}
and
\begin{equation}
    \frac{\partial q_{sat}}{ \partial b}(D^n,b^n) = \frac{\partial}{\partial b} \left( \frac{q_0 H}{D + B} e^{20(1-b/g)}\right)
    = -\frac{20}{g} q_{sat}(D^n, b^n).
\end{equation}
Then $\gamma_v$ is given by:
\begin{equation} \label{eq:gamma_v}
    \gamma_v = \frac{1}{1 + q_{sat}(D^n,b^n)\left[ \frac{20 \beta_2}{g} + \frac{\beta_1}{D + B}\right]}.
\end{equation}

Our choice for $\beta_1$ comes from Equation (8) in Vallis and Penn \cite{VP2020}. Taking $L = 2.4 \times 10^6 \text{ J kg}^{-1}$, $c_p = 1004 \text{ J kg}^{-1} \text{K}^{-1}$, $Q_a = 0.02$, $T_0 = 300 \text{K}$ and $H=10000\text{m}$ gives $\beta_1$ of approximately 1600 m. $\beta_2$ is $gL$ where $L \approxeq 10$ is the pseudo-latent heat from Zerroukat and Allen \cite{ZA2015}. Following Zerroukat and Allen \cite{ZA2015}, we take the threshold for rain $q_{precip}$ to be $1 \times 10^{-4}$ and the conversion rate for rain $\gamma_r$ to be $1 \times 10^{-3} \text{ s}^{-1}$. 

\section{Discretisation}

\subsection{Compatible finite element discretisation}
Our moist shallow water model is discretised in space using the compatible finite element method. Cotter and Shipton in \cite{CotterShipton2012} demonstrated that compatible finite element methods applied to the linear shallow water equations offer desirable properties analogous to the C-grid staggered finite difference method. These include the ability to support steady geostrophic modes, the avoidance of spurious mode branches, and an accurate representation of the dispersion relation for Rossby and inertia-gravity waves. We use an extension of this discretisation to the non-linear shallow water equations described in, for example, \cite{ShiptonCotterGibson} and \cite{Cotter2023}, and further extend it to include finite element spaces for buoyancy $b$ and moisture $q$ variables. The discretisation is based on a two-dimensional discrete de Rham complex $(\mathbb{V}^1, \mathbb{V}^2)$, where $\bm{u} \in \mathbb{V}^1$ and $D \in \mathbb{V}^2$. We use a mesh of triangular elements with the finite element family $\mathbb{V}^1 = BDM_2$ (Brezzi-Douglas-Marini elements on triangles) and $\mathbb{V}^2 = DG_1$, following the demonstration by Cotter and Shipton in \cite{CotterShipton2012} that these spaces can offer the desired properties described above. The choice means that $\bm{u}$ belongs in a space of quadratic vector-valued functions whose normal components are continuous across cell edges, and $D$ in a space of linear functions with no inter-element continuity constraints - the next-to-lowest order spaces for this finite element family. We have chosen to use the same discontinuous space for the buoyancy and moisture variables as we use for the depth, motivated by the fact that $D$ and $b$ appear in a similar way in the pressure gradient and the fact that using this space allows us to use the same transport scheme for these variables.

To write the weak forms for each equation we multiply by test functions from the relevant function spaces and integrate over the domain $\Omega$. We illustrate this for the velocity equation in the models with non-constant buoyancy, so that the right-hand side is $-b \nabla D -\frac{D}{2} \nabla b$. In the moist convective model the buoyancy is a constant $g$ and this term reduces to $-g \nabla D$, as in the usual shallow water equations. Thus the discretisation in that case reduces to that of the shallow water equations in, for example, \cite{ShiptonCotterGibson}. We write the velocity equation in vector invariant form, discretising:
\begin{equation}
    \frac{\partial{\bm{u}}}{\partial{t}} + (\nabla \times \bm{u}) \times \bm{u} + \frac{1}{2} \nabla \lvert \bm{u} \rvert^2 + f \bm{\hat{k}} \times \bm{u}= -b \bm{\nabla} (D + B) - \frac{D}{2} \bm{\nabla} b,
\end{equation}
by multiplying by a test function ${\bm{\psi}}$ from the $\mathbb{V}^1$ space and integrating. After integration by parts we have:
\begin{equation}
\begin{split}
    \int_{\Omega} {\bm{\psi}} \cdot \frac{\partial {\bm{u}}}{\partial t} \text{d}x
    + \int_{\Omega} \bm{u} \cdot \left(\nabla \times (\bm{u} \times \bm{\psi})\right) \text{d}x
    - \int_{\Gamma} \left( \bm{n} \times \llbracket \bm{u} \times \bm{\psi}\rrbracket \right) \cdot \tilde{\bm{u}} \text{d}S
    - \frac{1}{2}\int_{\Omega} \nabla \cdot \bm{\psi} \lvert \bm{u} \rvert^2 \text{d}x
    + \int_{\Omega} {\bm{\psi}} \cdot (f \bm{\hat{k}} \times \bm{u}) \text{d}x
    \\
    = \int_{\Omega} (B + D) \nabla \cdot ({\bm{\psi}}b) \text{d}x - \int_{\Gamma} \langle B + D\rangle \llbracket b{\bm{\psi}} \rrbracket \text{d}S
    + \frac{1}{2} \int_{\Omega} b \nabla \cdot (D {\bm{\psi}}) \text{d}x - \frac{1}{2} \int_{\Gamma} \langle b\rangle \llbracket D{\bm{\psi}} \rrbracket \text{d}S, \quad \forall \bm{\psi} \in \mathbb{V}^1, 
\end{split}
\end{equation}
where $\Gamma$ represents the set of cell edges, $\text{d}S$ is the measure for integrating over all cell edges, $\tilde{\bm{u}}$ is the upwind value of $\bm{u}$ and $\bm{n}$ is the unit normal pointing outward from the facet. Due to some fields having discontinuities across facets, their values on the facet are undefined. For the transported variable in the transport terms we choose the upwind value. For other terms we require the `jump' across the facet, denoted by $\llbracket \psi \rrbracket = \psi^+ - \psi^-$ (where each edge has sides arbitrarily labelled ``+" and ``-") and the average, denoted by $\langle \psi\rangle = (\psi^+ + \psi^-)/2$. This treatment of the pressure gradient term is similar to the treatment of the pressure gradient term in the discretisation of the compressible Euler equations of Melvin et al. \cite{GungHoCartesian}, Natale et al. \cite{Natale2016}, Bendall et al. \cite{Bendall2020compatible} and Cotter and Shipton \cite{CotterShipton2023}. The upwinded vector invariant form for the transport of $\bm{u}$ first appeared in \cite{NataleCotter2018} and was used for the rotating shallow water equations on the sphere in \cite{Gibson2020slate}. 

We write the weak form of the depth equation (Equation \eqref{eq:general_D_eq}) by multiplying by a test function $\phi$ from the $\mathbb{V}^2$ space and integrating, using integration by parts on the second term. Since $D$ is being transported, we define $D$ on cell edges as the upwind value with respect to the velocity, denoted by $\Tilde{D}$. This gives the weak form:

\begin{equation}
    \int_{\Omega} \phi \frac{\partial D}{\partial t} \text{d}x - \int_{\Omega} D {\bm{u}} \cdot  \nabla \phi  \text{d}x + \int_{\Gamma} \Tilde{D} \llbracket \phi \rrbracket {\bm{u}} \cdot \bm{n} \text{d}S = \int_\Omega \phi \beta_1 S_D \text{d}x, \quad \forall \phi \in \mathbb{V}^2,
\end{equation}
again where the `jump' terms are denoted by $\llbracket \cdot \rrbracket$.

In a similar way, the weak form of the buoyancy equation is:
\begin{equation}
    \int_{\Omega} \lambda \frac{\partial b}{\partial t} \text{d}x - \int_{\Omega} b \nabla \cdot (\lambda {\bm{u}}) \text{d}x + \int_{\Gamma} \Tilde{b} \llbracket \lambda {\bm{u}}\rrbracket \text{d}S = \int_{\Omega} \lambda \beta_2 S_b \text{d}x, \quad \forall \lambda \in \mathbb{V}^2,
\end{equation}
and the moisture equation is:
\begin{equation}
    \int_{\Omega} \tau \frac{\partial q^{(k)}}{\partial t} \text{d}x - \int_{\Omega} q^{(k)} \nabla \cdot (\tau {\bm{u}}) \text{d}x + \int_{\Gamma} \Tilde{q}^{(k)} \llbracket \tau {\bm{u}} \rrbracket \text{d}S = \int_\Omega \tau S_q^{(k)} \text{d}x, \quad \forall \tau \in \mathbb{V}^2.
\end{equation}

\subsection{Time stepper}
The prognostic variables are evolved in time using a semi-implicit quasi-Newton time-stepping scheme. This scheme follows the time stepper used by the Met Office's models, as described in \cite{Wood2014inherently} and \cite{GungHoCartesian}, and is also the time stepper used in \cite{Bendall2019recovered} and \cite{Bendall2020compatible}. The scheme separates the equation into terms that are treated explicitly and terms that are treated semi-implicitly. Each term in the equation is categorised as a forcing term $\mathcal{F}(\chi)$, a transport term $\mathcal{A}(\chi)$, or a physics term $\mathcal{P}(\chi)$, where $\chi = (\bm{u}, D, b, q^k)$ is a single vector encapsulating all prognostic variables. Forcing terms are treated semi-implicitly and transport and physics terms are treated explicitly, with physics happening separately to the dynamics at the end of the time-step.

Though the scheme has been described previously (see, for example, \cite{GungHoCartesian, Wood2014inherently, Bendall2019recovered, Bendall2020compatible}), below we discuss the linear solve in some detail, given that this part of the scheme involves a new thermal linear solver. In the description and the pseudocode, the prognostic variables that are evolved in time are denoted by the single vector $\chi$. The algorithm begins with a half timestep of explicit forcing, which takes the state at the $n$th time level, $\chi^n$, and returns $\chi_{f, e}$. This is given in weak form by:
\begin{equation}
\begin{split}
    \chi_{f, e} = \chi^n + \frac{1}{2}\Delta t \left[ -\int_{\Omega} {\bm{\psi}} \cdot (f \hat{\bm{k}} \times {\bm{u}}^n) \text{d}x
    +\int_{\Omega} (B + D^n) \nabla \cdot ({\bm{\psi}}b^n) \text{d}x 
    -\int_{\Gamma} \langle B + D^n\rangle \llbracket b^n{\bm{\psi}} \rrbracket \text{d}S \right.
    \\
    \left. +\frac{1}{2} \int_{\Omega} b^n \nabla \cdot (D^n {\bm{\psi}}) \text{d}x
    -\frac{1}{2} \int_{\Gamma} \langle b^n\rangle \llbracket D^n{\bm{\psi}} \rrbracket \text{d}S \right]
\end{split}
\end{equation}
and is written in the pseudocode as:
\begin{equation}
    \chi_{f, e} = \chi^n + \frac{1}{2}\Delta t \mathcal{F}(\chi^n),
\end{equation}
where $\mathcal{F}$ is the forcing operator which corresponds to the Coriolis and pressure gradient terms.

The next step is an outer iterative loop, which begins with computation of the advecting velocity $\bar{\bm{u}}$ as the average of the velocity at the current time level $\bm{u}^{n}$ and the best guess of the velocity at the next time level $\bm{u}^{n+1}$ : $\bar{\bm{u}} = \frac{1}{2} (\bm{u}^n + \bm{u}^{n+1})$. This advecting velocity is then used to explicitly advect the prognostic fields, $\mathcal{A}_{\bar{\bm{u}}}(\chi_{f, e})$, returning $\chi_T$. For the transport of $\bm{u}$, the equation is written in vector-invariant form and timestepped using the implicit midpoint method. All other fields are transported using an upwind discontinuous Galerkin scheme with a three-step Runge-Kutta timestepping procedure (SSPRK3), as outlined in, for example \cite{Bendall2020compatible} and \cite{Bendall2019recovered}. Within the Runge-Kutta step, we apply the limiter of Kuzmin \cite{Kuzmin2010vertex} as is also described by Cotter and Kuzmin \cite{CotterKuzminlimiters} to prevent the generation of new maxima and minima in the transport step.

This step is written in the pseudocode as:
\begin{equation}
    \chi_T = \mathcal{A}_{\bar{\bm{u}}}(\chi_{f, e}). \end{equation}
After that we start the inner loop, which begins with a half timestep of implicit forcing, given by:
\begin{equation}
    \begin{split}
    \chi_{f, i} = \chi_T + \frac{1}{2}\Delta t \left[ -\int_{\Omega} {\bm{\psi}} \cdot (f \hat{\bm{k}} \times {\bm{u}}^{n+1}) \text{d}x
    +\int_{\Omega} (B + D^{n+1}) \nabla \cdot ({\bm{\psi}}b^{n+1}) \text{d}x 
    -\int_{\Gamma} \langle B + D^{n+1}\rangle \llbracket b^{n+1}{\bm{\psi}} \rrbracket \text{d}S \right.
    \\
    \left. +\frac{1}{2} \int_{\Omega} b^{n+1} \nabla \cdot (D^{n+1} {\bm{\psi}}) \text{d}x
    -\frac{1}{2} \int_{\Gamma} \langle b^{n+1}\rangle \llbracket D^{n+1}{\bm{\psi}} \rrbracket \text{d}S \right]
\end{split}
\end{equation}
and in the pseudocode by:
\begin{equation}
    \chi_{f, i} = \chi_T + \frac{1}{2} \Delta t \mathcal{F}(\chi^{n+1}),
\end{equation}
where the $\mathcal{F}$ operator corresponds to the Coriolis and pressure gradient forces.
We define a residual $\chi_{r}$ as the difference between this newly-forced state and the best guess for the state at the next timestep $\chi^{n+1}$ by:
\begin{equation}
    \chi_r = \chi_{f, i} - \chi^{n + 1},
\end{equation}
and the aim then becomes reducing this residual to zero. A full Newton method would solve this iteratively as:
\begin{equation}
    \mathcal{J}\left[\chi^{n+1}_k\right]\left(\chi_{k+1}^{n+1} - \chi_{k}^{n+1}\right)  = - \chi_{r},
\end{equation}
where $\mathcal{J}$ is the Jacobian of $\chi_r$ with respect to $\chi^{n+1}$ and the subscript $k$ indicates the iteration index of the Newton loop. To avoid computing the Jacobian at every time step, we instead use a quasi-Newton method and approximate the Jacobian as a simple linear system, which no longer depends on $\chi_{k}^{n+1}$:
\begin{equation}
    \mathcal{J}\left[\chi_{k}^{n+1}\right] \approx \mathcal{S},
\end{equation}
so that we are instead solving
\begin{equation} \label{eq:linear_solve}
    \mathcal{S} \delta \chi = - \chi_{r},
\end{equation}
where $\delta \chi$ is the increment $\chi_{k+1}^{n+1} - \chi_{k}^{n+1}$.

Following Wood et al. \cite{Wood2014inherently} the choice of the linear operator $\mathcal{S}$ is inspired by the linearisation of $\chi_r$ about some reference state $\chi_{ref}$. The moisture variables are not included in the linearised equations. Excluding moisture from the linear solve matches the approach taken in the Met Office's ENDGame \cite{Wood2014inherently} and GungHo \cite{GungHoCartesian} models, as well as the compressible Euler model from Bendall et al. \cite{Bendall2020compatible}. The buoyancy variable is included in all formulations that include thermal effects. In these cases, we linearise the thermal shallow water equations rather than the shallow water equations. The strategy is to eliminate the buoyancy variable $b$ by writing it in terms of the velocity $\bm{u}$ and depth $D$, then to solve the resulting system for $\bm{u}$ and $D$ using a hybrid-mixed method before reconstructing $b$.

To derive our linear operator $\mathcal{S}$, we begin by linearising the thermal shallow water equations (which do not include physics terms),
\begin{align}
    &\frac{\partial{\bm{u}}}{\partial{t}} + (\bm{u} \cdot \nabla)\bm{u} + f \bm{\hat{k}} \times \bm{u} = -b \nabla (B + D) - \frac{D}{2} \nabla b,  \\ 
    &\frac{\partial{D}}{\partial{t}} + \nabla \cdot (\bm{u}D) = 0,  \\
    &\frac{\partial b}{\partial t} + (\bm{u} \cdot \nabla)b = 0,
\end{align}
about a background state $\chi = (\bm{u}_{ref}, H, b_{ref})$. $\bm{u}_{ref} =0$, $H$ is the mean background depth and $b_{ref}$ is taken as the initial buoyancy expression (so is constant in time but space-dependent), set when the model is initialised. Because the solve is approximate, we simplify the $\bm{u}$ equation by assuming no topography to give the linear system:
\begin{align}
    &\frac{\partial{\bm{u}}'}{\partial{t}} + f \bm{\hat{k}} \times \bm{u}' +b_{ref}\nabla D' + \frac{1}{2} \left[ H \nabla b_{ref} + H \nabla b' + D' \nabla b_{ref}\right] = 0, \\
    &\frac{\partial D'}{\partial t} + \nabla \cdot (\bm{u}' H) = 0, \\
    &\frac{\partial b'}{\partial t} + (\bm{u}' \cdot \nabla)b_{ref} = 0,
\end{align}
where the primed quantities represent perturbations from the background state.

Using the operator for this linear system to approximate the Jacobian, we obtain the matrix-vector problem corresponding to Equation \eqref{eq:linear_solve}:
\begin{align}
\begin{split}
    \int_\Omega {\bm{\psi}} \cdot \delta \bm{u} \text{d}x
    + \Delta t \int_{\Omega} {\bm{\psi}} \cdot (f \bm{\hat{k}} \times \delta \bm{u}) \text{d}x
    - \Delta t \int_\Omega  \nabla \cdot (\bar{b} {\bm{\psi}}) \delta D \text{d}x
    + \Delta t \int_{\Gamma} \llbracket \bm{\psi}\bar{b} \rrbracket \langle \delta D \rangle \text {d}S
    - \frac{\Delta t}{2}\int_\Omega H \bar{b} \nabla \cdot {\bm{\psi}} \text{d}x
    \\
    - \frac{\Delta t}{2}\int_\Omega H \delta b \nabla \cdot {\bm{\psi}} \text{d}x
    - \frac{\Delta t}{2}\int_\Omega \bar{b} \nabla \cdot (\delta D {\bm{\psi}}) \text{d}x
    + \frac{\Delta t}{2} \int_\Gamma \llbracket \delta D {\bm{\psi}} \rrbracket \langle \bar{b}\rangle \text{d}S = - \int_\Omega \bm{\psi} \cdot \bm{u}_r \text{d}x, \quad \forall \bm{\psi} \in \mathbb{V}^1,
    \\
    \int_\Omega \phi \delta D \text{d}x + \Delta t \int_\Omega \phi H \nabla \cdot \delta \bm{u} \text{d}x = - \int_\Omega \phi D_r \text{d}x, \quad \forall \phi \in \mathbb{V}^2,
    \\
    \int_{\Omega} \lambda \delta b \text{d}x - \Delta t \int_{\Omega} b_{ref} \nabla \cdot (\lambda {\delta \bm{u}}) \text{d}x + \Delta t \int_{\Gamma} \Tilde{b}_{ref} \llbracket \lambda {\delta \bm{u}}\rrbracket \text{d}S = - \int_\Omega \lambda b_r \text{d}x, \quad \forall \lambda \in \mathbb{V}^2,
\end{split}
\end{align}
where the increment we are solving for, $\delta \chi$, is given by $\delta \chi = (\delta \bm{u}, \delta D, \delta b)$. We solve this system by eliminating the buoyancy equation via:
\begin{equation} \label{eq:b_prime}
    \delta b = -\Delta t (\delta \bm{u} \cdot \nabla) b_{ref} + b_r,
\end{equation}
to give a reduced system for $\bm{u}$ and $D$.
This reduced system is solved using a hybridisation method, which is described in detail in \cite{Firedrakebook, Gibson2020slate} as well as in \cite{Bendall2020compatible}. The technique involves relaxing the continuity requirements on the normal components of functions in $\mathbb{V}^1$ and enforcing continuity through Lagrange multipliers as part of the solution formulation. We then back-substitute using these $(\delta \bm{u}, \delta D)$ values to get a value for $\delta b$. In the moist convective formulation there is no prognostic buoyancy field, so our operator $\mathcal{S}$ is derived by linearising the standard shallow water equations, resulting in a system for $\delta \bm{u}$ and $\delta D$ which is solved in the same way as the reduced system above. These increments $\delta \chi = (\delta \bm{u}, \delta D, \delta b)$ (or $\delta \chi = (\delta \bm{u}, \delta D)$ in the moist convective model) are added to the current guess to update $\chi^{n + 1}$ and the iteration begins again. 

In practice we do not iterate to convergence, but rather perform a fixed number of iterations of both this inner loop and the outer loop. After exiting both loops the physics processes are treated explicitly:
\begin{equation}
    \chi^{n+1} = -\int_\Omega \phi \beta_1 P \text{d}x
    + \int_\Omega \lambda \beta_2 P \text{d}x
    + \int_\Omega \tau S_q \text{d}x
\end{equation}
written in the pseudocode as:
\begin{equation}
    \chi^{n+1} = \mathcal{P}(\chi^{n+1})
\end{equation}
In the last step of the algorithm the timestep is advanced:
\begin{equation}
    \chi^n = \chi^{n+1}
\end{equation}

\begin{algorithm}
	\caption{Semi-implicit quasi-Newton time-stepping scheme} 
	\begin{algorithmic}[1]
		\State Set: $\chi^{n+1} = \chi^{n}$
        \State Explicit forcing: $\chi_{f, e} = \chi^n + \frac{1}{2}\Delta t \mathcal{F}(\chi^n)$
			\For {$OUTER$}
				\State Update advecting velocity: $\bar{\bm{u}} = \frac{1}{2} (\bm{u}^{n+1} + \bm{u}^n)$
				\State Explicit transport: $\chi_T = \mathcal{A}_{\bar{\bm{u}}}(\chi_{f, e})$
                \For {$INNER$}
                    \State Implicit forcing: $\chi_{f, i} = \chi_T + \frac{1}{2}\Delta t \mathcal{F}(\chi^{n+1})$
                    \State Write residual: $\chi_{r} = \chi_{f, i} - \chi^{n+1}$
                    \State Solve: $\mathcal{S}(\delta \chi) = \chi_{r}$ for $ \delta \chi$
                    \State Increment: $\chi^{n+1} = \chi^{n+1} + \delta \chi$
			    \EndFor
            \EndFor
			\State Explicit physics: $\chi^{n+1} = \mathcal{P}(\chi^{n+1})$
			\State Advance time step: $\chi^n = \chi^{n+1}$
	\end{algorithmic} 
\end{algorithm}

\section{Results}

In this section we give details of three test cases in each moist shallow water formulation, along with results. Our model is written using Gusto, the dynamical core toolkit built on the Firedrake finite element library. Firedrake provides automated code generation for the solution of partial differential equations using the finite element method \cite{FiredrakeUserManual}. Code is automatically generated using the Unified Form Language \cite{UFL} and equations are provided to PETSc \cite{PETSc}, which provides direct access to runtime configurable iterative solvers and preconditioners. The hybridisation of the linear system in the linear solve stage is implemented using the Slate framework \cite{Gibson2020slate}, which performs local elimination and recovery options. The reduced equation for the trace variables is numerically inverted using the conjugate gradient method and PETSc's GAMG preconditioner.

Our first step to verify our model was to run Test 2 from \cite{ZA2015} using the moist thermal formulation, with the initial conditions as given there. We found that it was necessary to reproduce as closely as we could the resolution and timestep in \cite{ZA2015} in order to see the same amount of cloud and rain. We also found that the amount of cloud and rain were both sensitive too to the use of the limiters in the transport step, and that more of both species were generated when no limiters were used. This experiment, while establishing confidence in our model, also demonstrated the sensitivity of the cloud field to aspects of the model and various modelling choices. (These results are not presented here.)

For the Earth-related physical constants for all tests we take rotation rate of the Earth $\Omega = 7.292 \times 10^{-5} \text{s}^{-1}$ where $f = 2 \Omega \cos{\phi}$ with $\phi$ the latitude, acceleration due to gravity $g = 9.80616 \text{ m } \text{s}^{-2}$ and radius of the Earth $R = 6371220 \text{ m}$.
As described in Zerroukat and Allen \cite{ZA2015}, in all our tests the initial vapour is dependent on the saturation, and is chosen to be near saturation so that the physics is activated easily and cloud and rain are produced without having to run the model for a long time. The scaling factor $q_0$ that appears in the saturation function $q_{sat}$ is chosen so that the maximum value of the initial vapour is around the atmospheric value of 0.02. For both the steady state test and the flow over a mountain test that means a value for $q_0$ of 0.007, and a value of 0.0027 for the unstable jet test. Rain, when it is produced, is not advected with the flow but instead accumulates where it forms.

For the moist convective tests the saturation function varies only with depth, as given by \eqref{eq:sat_func_depth}, but includes a time-constant, spatially-varying factor $\theta$. $\theta$ has a similar form to the potential temperature initial condition in Test 1 of Zerroukat and Allen \cite{ZA2015}  \footnote{Note the sign change in the numerator relative to the equation given in Zerroukat and Allen \cite{ZA2015}, which was necessary to achieve balanced initial conditions.}:
\begin{equation} \label{eq:theta}
    \theta = \frac{\theta_0 + \sigma \text{cos}^2\phi \left[ (\omega + \sigma) \text{cos}^2\phi + 2(\Phi_0 - \omega - \sigma)\right]}{\Phi_0^2 + (\omega + \sigma)^2 \text{sin}^4 \phi -2 \Phi_0 (\omega + \sigma)\text{sin}^2 \phi},
\end{equation}
with $u_0 = 20 \text{ m s}^{-1}$, $\omega = \left( \Omega a u_0 + \frac{u_0^2}{2}\right)$, $\sigma = \omega/10$, $\Phi_0 = 3 \times 10^4 \text{ m}^2 \text{ s}^{-2}$, and $\theta_0 = \epsilon \Phi_0^2$ with $\epsilon = 1/300$.

All tests are run on an icosahedral grid. The moist flow over a mountain test uses a grid of 20480 cells with a maximum cell edge length of 263 km and a minimum cell edge length of 171 km. The moist unstable jet test uses a grid of 81920 cells with a maximum cell edge length of 132 km and a minimum cell edge length of 85 km.

\subsection{Moist Steady State Geostrophic Flow}
The first test is a modification of test 2 of Williamson et al. \cite{Williamson1992}, and is based on the steady state test of Zerroukat and Allen \cite{ZA2015} (written there in terms of geopotential rather than depth, and potential temperature rather than buoyancy)\footnote{Note the sign change in Equation \eqref{eq:W2_buoyancy} relative to the potential temperature equation in Zerroukat and Allen \cite{ZA2015} to ensure balanced initial conditions.}. The test consists of zonally balanced flow on a sphere with no topography. The initial conditions are solutions to the thermal shallow water equations and so are a slight modification of the original Williamson et al. \cite{Williamson1992} initial conditions to account for the extra terms in the thermal system. Because the test is initialised with a steady state, the expectation is that the solution should remain at the initial state throughout the simulation, and thus the metric for evaluating the results is how well the final state compares to the initial state. The model is initialised with water vapour $q_v$ either at or close to saturation everywhere and no cloud or rain. The initial conditions for both of the models with buoyancy (moist thermal and moist convective thermal) are (with $(\lambda, \phi)$ representing spherical (longitude, latitude) coordinates):
\begin{align}
    &\bm{u}(\lambda, \phi) = (u_0 \text{cos} \phi, 0), \\
    &D(\lambda, \phi) = H - \frac{1}{g} (\omega + \sigma) \text{sin}^2\phi, \\
    &b(\lambda, \phi) = g \left( 1 - \frac{\theta_0 + \sigma \text{cos}^2\phi \left[ (\omega + \sigma) \text{cos}^2\phi + 2(\Phi_0 - \omega - \sigma)\right]}{\Phi_0^2 + (\omega + \sigma)^2 \text{sin}^4 \phi - 2 \Phi_0 (\omega + \sigma) \text{sin}^2 \phi}\right),
    \label{eq:W2_buoyancy}\\
    &q_v(\lambda, \phi) = (1 - \xi) q_{sat}(D, b),
\end{align}
where $u_0 = 20 \text{ m s}^{-1}$, $\omega = \left( \Omega R u_0 + \frac{u_0^2}{2}\right)$, $\sigma = \omega/10$, $\Phi_0 = 3 \times 10^4 \text{ m}^2 \text{ s}^{-2}$, and $\theta_0 = \epsilon \Phi_0^2$ with $\epsilon = 1/300$. The background depth $H$ is $\Phi_0/g$. $\xi$ is a parameter that controls how near to saturation the initial vapour is. When $\xi = 0$ the vapour is initialised with the initial saturation values, and when $\xi = 10^{-3}$ vapour is initialised just below saturation. 

In the moist convective test we use the same initial condition for $\bm{u}$ and we take the expression for $D$ from the thermal case but with $\sigma=0$. Setting $\sigma=0$ in the $D$ expression excludes the terms required to balance the buoyancy equation and reduces the initial depth to that of the original Williamson et al. test \cite{Williamson1992}, meaning balanced initial conditions in this set-up where there is no buoyancy field. The saturation function varies only with depth and is given by Equation \eqref{eq:sat_func_depth}. Note that the $\theta$ used in this saturation function expression (given in Equation \eqref{eq:theta}) has a non-zero $\sigma$, which ensures that the initial vapour has a similar form across all versions of the test.

\begin{figure}[htp]
    \centering
    \includegraphics[width=\textwidth]{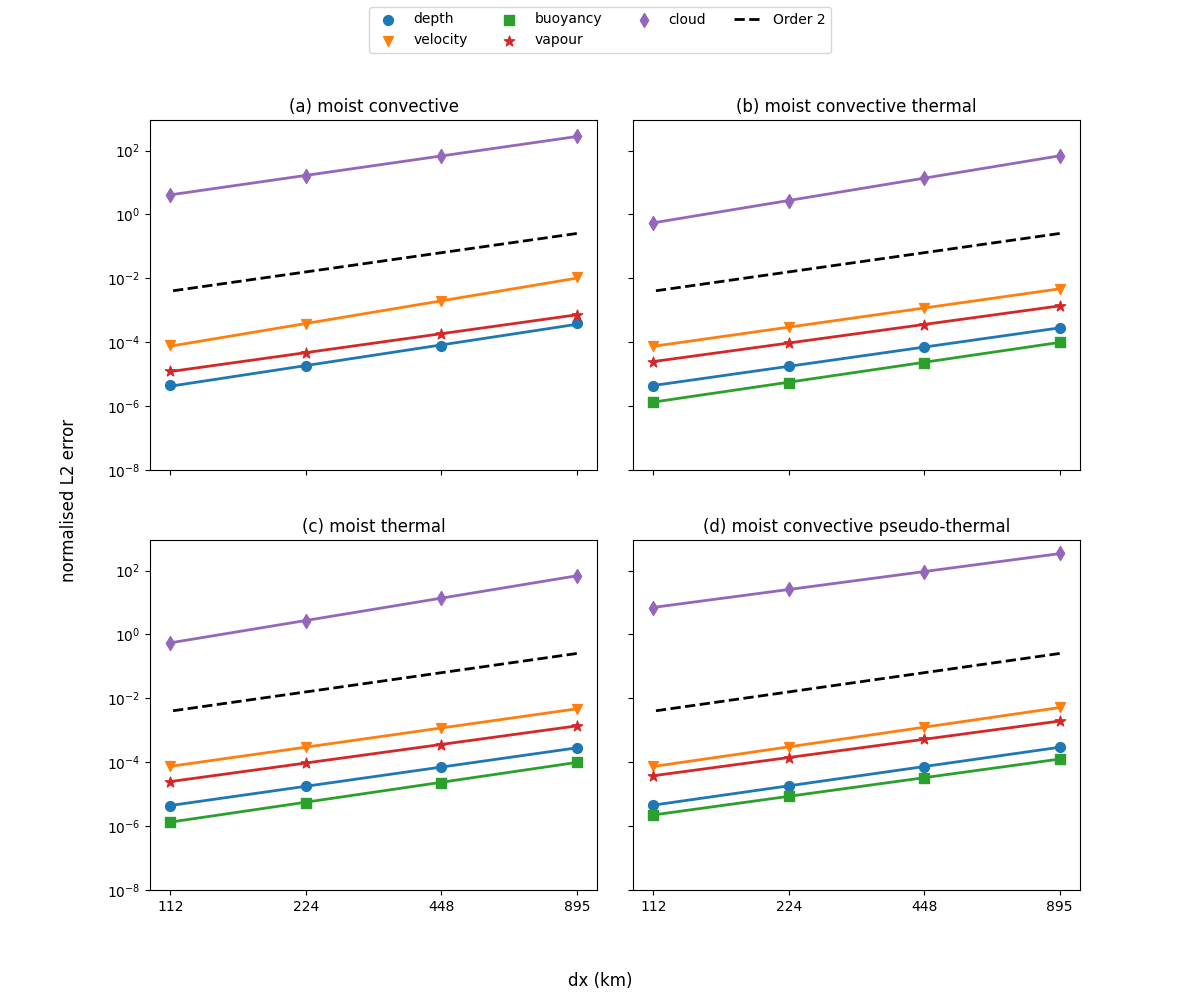}
    \caption{Convergence plots for the moist Williamson 2 tests in each framework. (a): the moist convective formulation, (b): the moist convective thermal formulation, (c): the moist thermal formulation, and (d): the moist convective pseudo-thermal formulation. The $L^2$ norm of the error is normalised by dividing it by the $L^2$ norm of the field in each case except for the cloud field, which is not normalised (because the test begins with no cloud). Second-order convergence lines are also shown on the plot for comparison.}
    \label{fig:W2_convergence}
\end{figure}

Because this is a steady state simulation, we show results from convergence tests for this example in each of the four frameworks. Following Williamson et al. \cite{Williamson1992}, we compute $L^2$ error norms for each field after 5 days at varying spatial resolutions. We vary the timestep at each resolution so as to give a constant advective Courant number of approximately 0.02, chosen to be small so as to generate accurate reference solutions. We choose $\xi = 0$ so that the model is initialised with vapour everywhere at saturation. This means that any numerical errors that produce overshoots will push vapour above saturation and generate clouds. The idea is that these numerical errors are sufficiently small that the amount of cloud produced stays below the threshold to produce rain and does not affect the stability of the overall flow. These convergence plots are shown in Figure \ref{fig:W2_convergence}. Because of the polynomial degree used in the finite elements, we expect second-order convergence from our time-stepping scheme, and second-order convergence lines are shown on all plots. All fields in each of the four frameworks converge at the expected rate.

\subsection{Moist Flow over an Isolated Mountain}
The second test is a modification of the flow over a mountain test of Williamson et al. \cite{Williamson1992} and is based on the version presented in Zerroukat and Allen \cite{ZA2015}. Unlike the steady state test, this test does not have an analytical solution. The purpose of the test is to assess how well the model can capture interactions between the moist physics and the dynamics, and how well cloud and rain are produced in a physically realistic way (for example, when warm air moves into a colder region). The initial velocity and depth are as in the original Williamson 5 test \cite{Williamson1992}, with the addition of $\sigma$ in the depth equation in all thermal models. The thermal versions of the test use the same expression for initial buoyancy as the steady state test, in contrast to Zerroukat and Allen \cite{ZA2015} where a new initial temperature was used. We made this choice because the buoyancy expression in the steady state test is derived from balances between velocity, depth and buoyancy. Including $\sigma$ in the depth equation in this test makes the depth and velocity fields the same as the steady state test (excluding topography) and hence we can use the buoyancy from that test to again give balanced initial conditions.
The topography is given by:
\begin{equation}
    B = h_0\left(1 - \frac{1}{R} \text{min} \left[ R, \sqrt{(\lambda - \lambda_c)^2 + (\phi - \phi_c)^2} \right]\right),
\end{equation}
with $h_0 = 2000$ m, $R = \pi/9$, $\lambda_c = 3\pi/2$ and $\phi_c = \pi/6$.
The initial conditions for all thermal models are:

\begin{align}
    &\bm{u}(\lambda, \phi) = (u_0 \text{cos} \phi, 0), \\
    &D(\lambda, \phi) = H - \frac{1}{g} (\omega + \sigma) \text{sin}^2\phi - B, \\
    &b(\lambda, \phi) = g \left( 1 - \frac{\theta_0 + \sigma \text{cos}^2\phi \left[ (\omega + \sigma) \text{cos}^2\phi + 2(\Phi_0 - \omega - \sigma)\right]}{\Phi_0^2 + (\omega + \sigma)^2 \text{sin}^4 \phi - 2 \Phi_0 (\omega + \sigma) \text{sin}^2 \phi}\right), \\
    &q_v = (1-\xi) q_{sat}(D, b),
\end{align}
with all constants as given in the steady state test except for $\xi$, which is taken as 0.02 here so that vapour is initialised close to saturation. The background depth $H$ is 5960 m and the test is run for 50 days.

Like the steady state test, in the moist convective version of this test the initial condition for $\bm{u}$ is the same as in the thermal case, the initial condition for $D$ is the same but with $\sigma = 0$, and the initial vapour is also $(1 - \xi)q_{sat}$ where $q_{sat}$ is now just a function of depth.

\begin{figure}[htp]
    \centering
    \includegraphics[width=\textwidth]{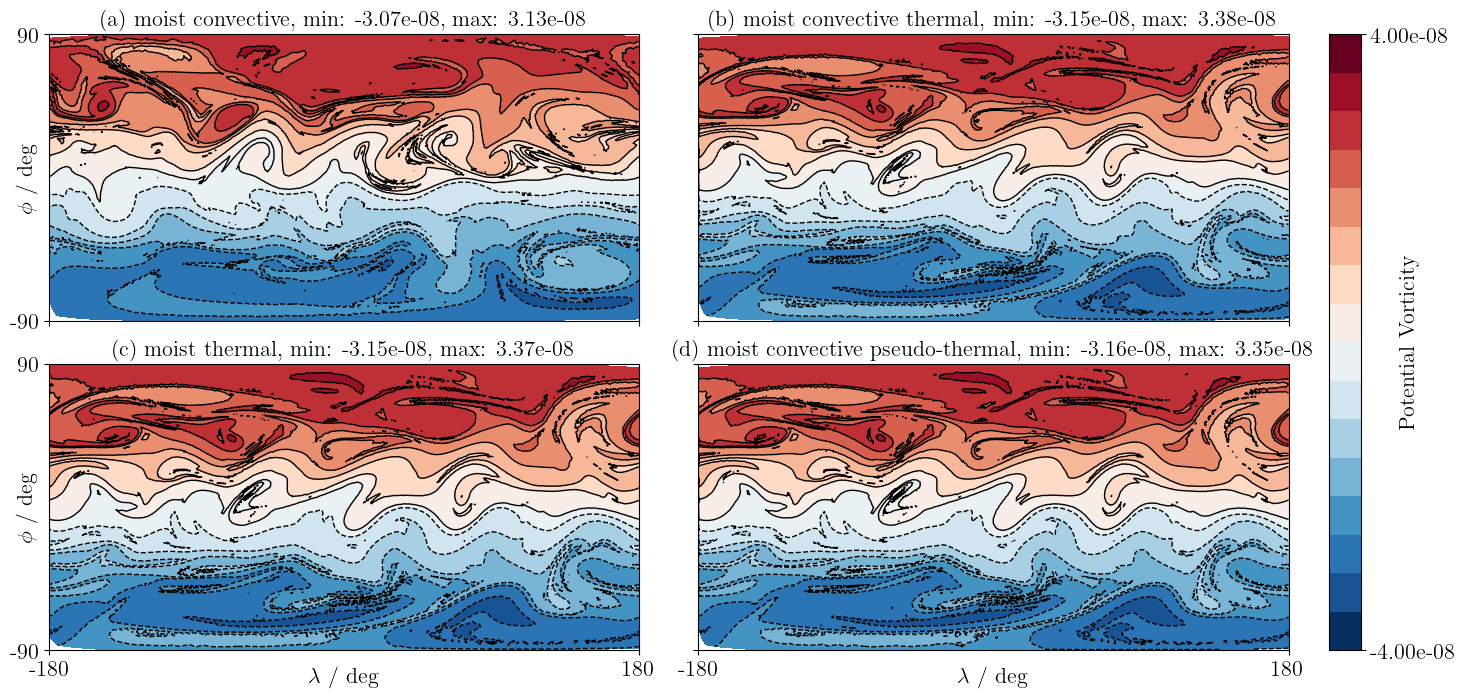}
    \caption{Potential vorticity (defined by $\frac{(\nabla \times \bm{u}) \cdot \bm{\hat{k}} + f}{D}$) at Day 50 in the flow over a mountain test. (a) the moist convective formulation, (b) the moist convective thermal formulation, (c) the moist thermal formulation, and (d) the moist convective pseudo-thermal formulation.}
    \label{fig:W5_vorticity}
\end{figure}

Figure \ref{fig:W5_vorticity} shows the potential vorticity after 50 days in the four different moist shallow water formulations. The large-scale flow dynamics are broadly similar across all the formulations, generated by similar flow patterns as the air is forced up and over the mountain. The formulation that is the most different to the others is the moist convective formulation (subplot a), which is the only formulation without a prognostic buoyancy field. The difference between this formulation and the moist convective pseudo-thermal model (subplot d) demonstrates that the presence of the buoyancy term in the $\bm{u}$ equation has a significant impact. Comparing the moist convective thermal formulation (subplot b) to the moist thermal formulation (subplot c) allows us to isolate the effect of the convective forcing on the depth equation, with the size of this effect being controlled by the size of the $\beta_1$ parameter. For our choice of $\beta_1$ these two formulations are very similar. A discussion of the effect of varying the size of the $\beta_1$ parameter is given in Appendix B. Comparing the moist convective thermal formulation (subplot b) to the moist convective pseudo-thermal formulation (subplot d) allows us to isolate the effect of the latent heat feedback on the buoyancy equation, as the only difference between these two models is the fact that the moist convective pseudo-thermal model sees no impact on the buoyancy equation from latent heat release and the moist convective thermal model does. Again we see that this effect is relatively small with our choice of the $\beta_2$ parameter.

\begin{figure}[htp]
    \centering
    \includegraphics[width=\textwidth]{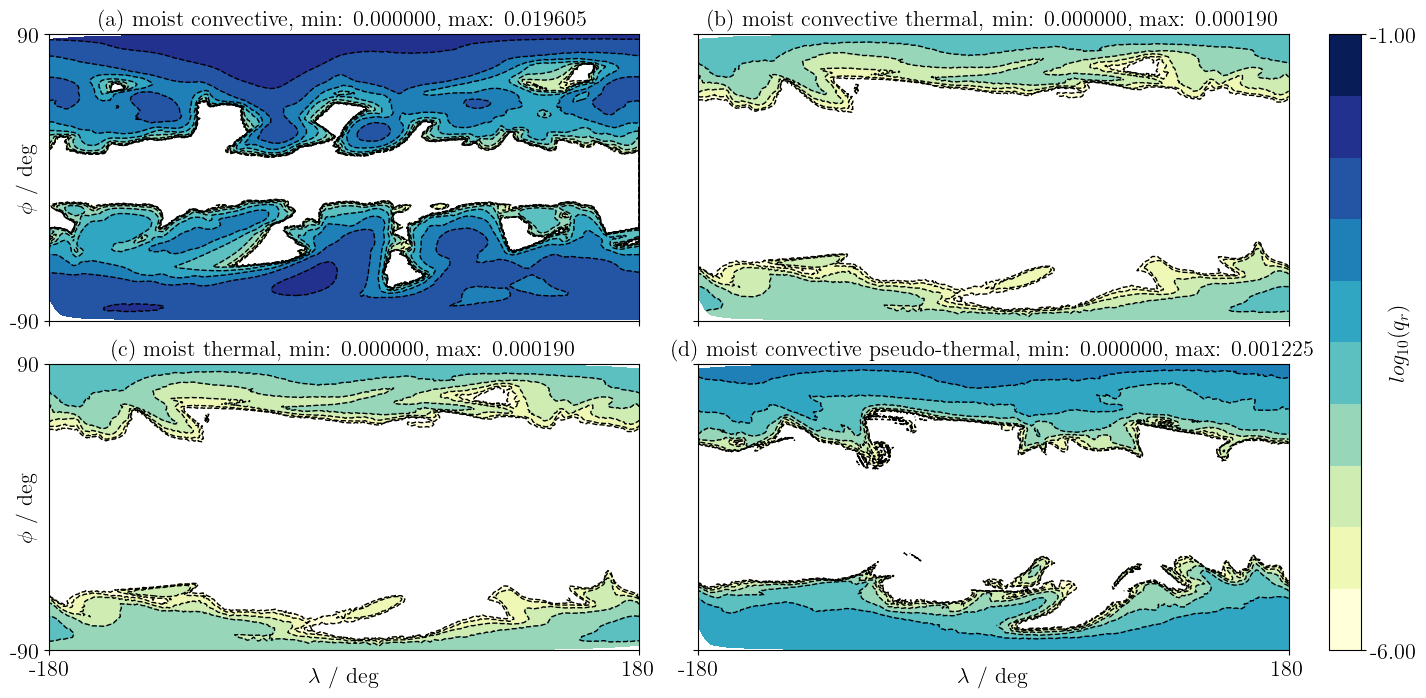}
    \caption{Accumulated rain at Day 50 in the flow over a mountain test. (a) the moist convective formulation, (b) the moist convective thermal formulation, (c) the moist thermal formulation, and (d) the moist convective pseudo-thermal formulation. We plot $\text{log}_{10}(q_r)$ rather than $q_r$ directly so that we can plot all four configurations on the same scale despite variation in rain quantities, and omit rain values of less than $1 \times 10^{-6}$ from the plot.}
    \label{fig:W5_rain}
\end{figure}

Figure \ref{fig:W5_rain} shows the accumulated rain after 50 days in each formulation. The differences between formulations seen in the potential vorticity results are even clearer here, and the differences in rain production mechanisms are also evident. Across the formulations cloud and rain are produced when warm air is blown into colder regions, and this pattern of accumulated rainfall at higher latitudes is evident in each formulation. The formulations with the most rain are the moist convective (subplot a) and the moist convective pseudo-thermal (subplot d). This can be explained by the fact that neither of these models see any heating effect from a moisture phase change. The latent heat feedback on the buoyancy in the other two models adjusts the buoyancy and hence the saturation function such that there is less cloud formed at the following timestep. Thus the latent heat feedback has a heating effect which acts to suppress further cloud formation, and the models with no latent heat effect go on to produce more cloud and later rain. There is more rain again in the moist convective model (subplot a), which we suggest is due to the fact that this formulation does not have a prognostic buoyancy variable. In all other formulations the buoyancy variable obeys the advective transport equation, meaning that the buoyancy and the vapour (which also obeys the advective transport equation) of a parcel of air will generally stay constant as the air parcel is blown around. In the moist convective model, however, the buoyancy-like variable is not a prognostic field and varies only in latitude. This background field appears in the saturation function and so as air parcels are blown meridionally they encounter a new saturation value which could cause vapour to suddenly exceed saturation and be converted to cloud and eventually rain.

\subsection{Moist Unstable Jet}
The third test is a modification of the barotropic instability test case of Galewsky et al. \cite{Galewsky2004} and follows the moist version of this test by Ferguson, Jablonowski and Johansen \cite{FergusonJablonowski2019}. In the dry version the test is initialised with a balanced zonal jet to which a small height perturbation is added. The height perturbation initiates the roll-up of the jet into well-defined vortices after a few days. In the moist version of the test these roll-ups trigger a moist response. The initial velocity is defined in Galewsky et al. \cite{Galewsky2004} and the initial depth is obtained by numerically integrating the balance equation, using this velocity (Equation (3) in that paper). The base flow is then perturbed to initiate the roll-up of the jet by adding a localised bump to the balanced depth field. The introduction of the buoyancy variable in the thermal model modifies the balance equation and means that the depth computed in Galewsky et al. \cite{Galewsky2004} no longer balances the velocity. To begin with a balanced jet we therefore solve the thermal balance equation for the depth $D$:

\begin{equation}
    u^2 \tan(\phi) \frac{1}{R} + f u = -\frac{1}{R} \left(-b \frac{\partial D}{\partial \phi} - \frac{D}{2} \frac{\partial b}{\partial \phi} \right),
\end{equation}
choosing $u$ as given in Galewsky et al. \cite{Galewsky2004} and defining the initial buoyancy as $b = g - \Delta b \cos(\phi)$ where $\phi$ is the latitude and $\Delta b$ is a constant. Solving this balance equation involves taking an integrating factor $G = \left (g - \Delta b \cos(\phi) \right )^{\frac{1}{2}}$ and then numerically integrating the equation:
\begin{equation}
    D = (g - \Delta b \cos(\phi))^{-\frac{1}{2}}\left[H (g-1)^{\frac{1}{2}} - \int (g - \Delta b \cos(\phi))^{-\frac{1}{2}} (u^2 \tan(\phi) + Rfu) \text{d}\phi \right],
\end{equation}
where $H (g-1)^{\frac{1}{2}}$ is an integration constant, chosen so that $D$ has the correct dimensions, and we take $\Delta b = 1 \text{ m s}^{-2}$. The balanced initial state is perturbed by adding a localised bump to the depth field, as in Galewsky et al. \cite{Galewsky2004}. The initial vapour field is set just below saturation:
\begin{equation}
    q_v(\lambda, \phi) = (1 - \xi) q_{sat}(D, b),
\end{equation}
with $\xi=0.02$. The initial cloud and rain fields are set to zero. The background depth $H$ is 10,000 m and we run the test for 6 days. 

There is no prognostic buoyancy variable in the moist convective formulation of this test so the initial velocity and depth in that formulation are the same as in the original jet in Galewsky et al. \cite{Galewsky2004}. The $\theta$ that appears in the saturation function \eqref{eq:sat_func_depth} is given by $\theta = -\frac{\Delta b \cos \phi}{g}$ (coming from the relationship between $b$ and $\theta$, with $\theta$ as given for the thermal versions of the test and $\Delta b$ again a constant $1 \text{ m s}^{-2}$). The initial vapour is, as before, $(1 - \xi) q_{sat}$, with $\xi = 0.02$. 

Like the flow over a mountain test there is no analytical solution for this test, but rather the results are assessed on the basis of the appearance of expected physical features, such as rain being generated along the fronts as the jet rolls up. The development and evolution of the jet instability in all four frameworks is consistent with the dry simulation results in Galewsky et al. \cite{Galewsky2004}, and with the moist thermal results in \cite{FergusonJablonowski2019}. In all four formulations we see cloud formation about 4.5 days into the simulation as cloud is produced in the vortices. By about 5.25 days the conversion from cloud to rain has begun, and by day 6 there is a clear pattern of cloud and rain matching the roll-ups of the jet.

\begin{figure}[htp]
    \centering
    \includegraphics[width=\textwidth]{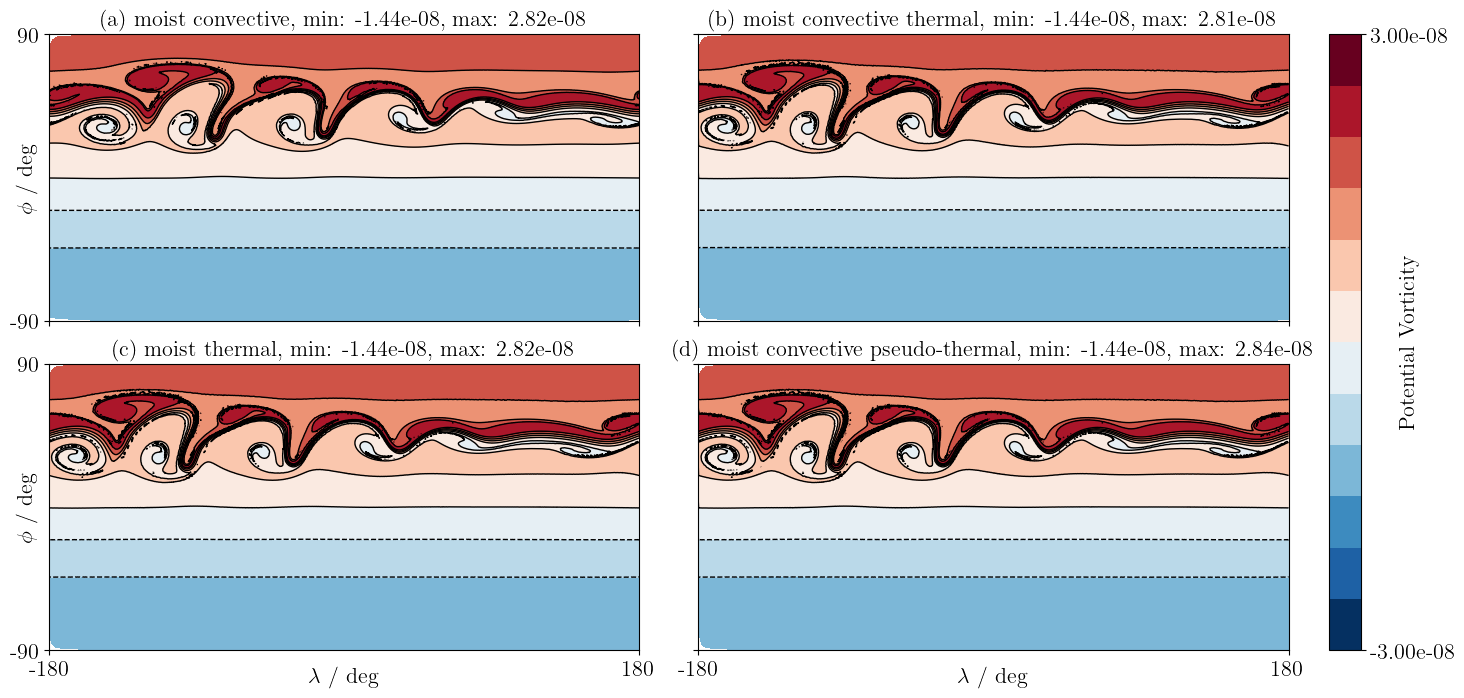}
    \caption{Potential vorticity (defined by $\frac{(\nabla \times \bm{u}) \cdot \bm{\hat{k}} + f}{D}$) at Day 6 in the moist unstable jet test. (a) the moist convective formulation,
    (b) the moist convective thermal formulation, (c) the moist thermal formulation, and (d) the pseudo-moist
    convective thermal formulation.}
    \label{fig:compare_vorticity_GJ}
\end{figure}

Figure \ref{fig:compare_vorticity_GJ} shows the potential vorticity after 6 days in each of the different moist shallow water formulations. The potential vorticity looks quite similar across the formulations, with the clearest difference in the moist convective formulation (subplot a), which arises because of the slightly different initial depth field used in this formulation. 

\begin{figure}[htp]
    \centering
    \includegraphics[width=\textwidth]{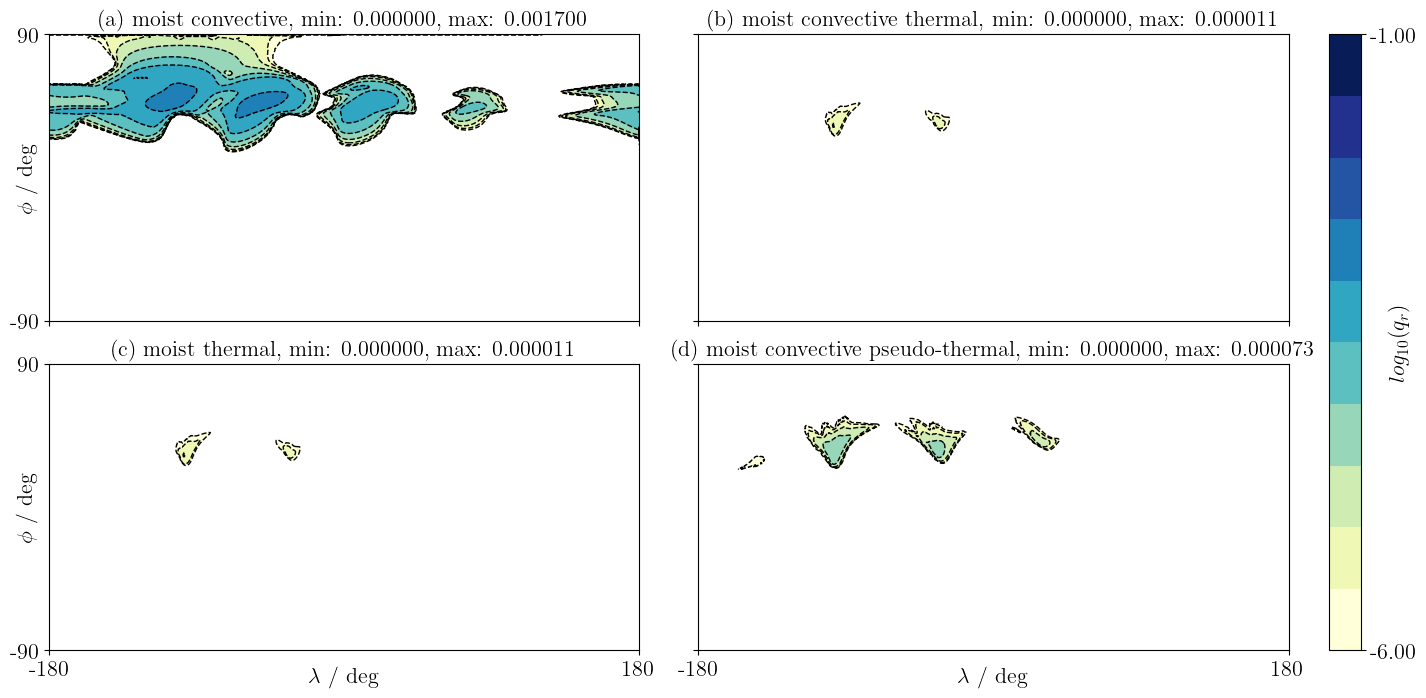}
    \caption{Accumulated rainfall at Day 6 in the moist unstable jet test. (a) the moist convective formulation,
    (b) the moist convective thermal formulation, (c) the moist thermal formulation, and (d) the pseudo-moist
    convective thermal formulation. As in Figure \ref{fig:W5_rain}, we plot $\text{log}_{10}(q_r)$ rather than $q_r$ directly so that we can plot all four configurations on the same scale despite variation in rain quantities, and omit rain values of less than $1 \times 10^{-6}$ from the plot.}
    \label{fig:compare_rain_GJ}
\end{figure}

Figure \ref{fig:compare_rain_GJ} shows the accumulated rainfall in this test case after 6 days in each formulation. As in the moist flow over a mountain test, the rain field is more discriminating than the potential vorticity field. In all formulations we see rain production in the roll-ups of the jet but the amount of rain varies between formulations. The rain suppression from buoyancy feedback that we saw in the mountain test is also visible here, where we again see the smallest rainfall amounts in both formulations with latent heat feedback effects (the moist convective thermal formulation (subplot b) and the moist thermal formulation (subplot c)). The moist convective formulation (subplot a) again has the most rain, which we attribute as before to the fact that the temperature-like variable in this formulation is constant in time. This is in contrast to the prognostic buoyancy variable in all of the other formulations, which makes it less likely for a parcel of air to suddenly exceed the saturation value and trigger condensation.

\section{Conclusions}
In this paper we have presented the first compatible finite element solutions to the moist shallow water equations, which is a significant step towards advancing the use of the moist shallow water model in physics-dynamics coupling investigations. The shallow water dynamics are computationally inexpensive and coupling them to the moist physics scheme offers a simple yet reasonably realistic model that includes such numerical complexities as non-linear switches and thresholds and new timescales for moist physics processes, which are highly relevant for physics-dynamics coupling strategies. We discretise a general formulation of the moist shallow water model, which allows for the representation of moisture effects as thermal, as convective, or as both thermal and convective. This means that we can recover three existing moist shallow water models with our implementation: the moist convective model of Bouchut et al. \cite{Bouchut2009}, the moist thermal model of Zerroukat and Allen \cite{ZA2015} or the moist convective thermal model of Kurganov, Liu and Zeitlin \cite{KLZ2020}. Our model also facilitates a new, fourth formulation of the moist shallow water equations, previously unexplored, which we call the moist convective pseudo-thermal shallow water model. We couple the moist shallow water model to a three-state moist physics scheme, where reversible conversions between water vapour and cloud feed back on the dynamics via a heating effect, a mechanical effect, or both, depending on the framework. 

Our new compatible finite element discretisation extends the discretisation for the shallow water equations from Cotter and Shipton \cite{CotterShipton2012} to include spaces for buoyancy and moisture variables. The compatible finite element method offers a promising solution to the problems of communication bottlenecks associated with the convergence of grid points at the poles in a latitude-longitude grid. This is particularly relevant as HPC systems move towards massively parallel set-ups and is prompting the move by the Met Office to compatible finite elements for their next-generation dynamical core. We have demonstrated the effectiveness of our discretisation by applying the model in all four frameworks to three different test cases: moist steady state geostrophic flow, moist flow over an isolated mountain, and a moist unstable jet. We have detailed the set-ups for all of these test cases in each formulation and given reference solutions. The purpose of this is two-fold: firstly, to demonstrate that our flexible model discretised using compatible finite elements gives realistic behaviour and cloud and rain formation, and secondly, to provide reference solutions with comparisons between the different moist shallow water models that convey information on the implications of different modelling choices. Our tests demonstrate the effect on the large-scale dynamics of both convective forcing on the depth field and latent heat feedback on the buoyancy field. The latent heat feedback on buoyancy has a significant impact on the amount of rain produced - in models with latent heat feedback we see that this feedback acts to suppress cloud and rain production. It is also evident that models with a prognostic temperature-like field see less cloud and rain formation than models with a time-constant background temperature field, where it is more likely for a parcel of air to suddenly exceed the saturation value because it does not have a conserved temperature-like variable. The results highlight the usefulness of our new moist convective pseudo-thermal model as a connection between the other three formulations, allowing us as it does to isolate the effects of modelling choices in this way.

The moist shallow water model provides a useful tool with which to investigate questions about physics-dynamics coupling in time-stepping schemes. Our compatible finite element discretisation of the general moist shallow water model offers a framework within which four different moist shallow water models can be interrogated and compared, including a new, previously unexplored moist convective pseudo-thermal shallow water model. We intend to follow this work up by using the moist shallow water equations to investigate how time steppers compare at coupling different processes in time. In addition to this, the test set-ups and example solutions we have detailed here can be useful as solutions to verify against for others working with moist shallow water models, irrespective of the numerical schemes employed.

\section*{Funding Acknowledgements}
N.H was supported by a Natural Environment Research Council (NERC) GW4+ Doctoral Training Partnership 2 Studentship (NE/S007504/1). For the purposes of open access, the authors have applied a creative commons attribution (CC BY) licence to any author accepted manuscript version arising from this submission.

\section*{Data Availability}
The data that support the findings of this study are available from the corresponding author upon reasonable request.

\section*{Conflicts of Interest}
The authors have no conflict of interest to declare.

\appendix
\section{Physics schemes}
Though all our tests are run using the three-state physics scheme from Zerroukat and Allen \cite{ZA2015}, we also implemented the one-way physics scheme described by Bouchut et al. \cite{Bouchut2009}. This is the scheme previously used with the moist convective and the moist convective thermal shallow water equations. To compare the two schemes we ran the moist convective version of the flow over a mountain test using each scheme in turn, defining the same depth-dependent saturation function in each case. In the one-way physics scheme this saturation function is the threshold for the direct conversion from vapour to rain and is not reversible. The proportion above saturation converted, $\gamma_r$, is set at $1 \times 10^{-3} \text{ s}^{-1}$. In the three-state scheme the same saturation function dictates conversion between vapour and cloud, which is a reversible conversion. The proportion converted, $\gamma_v$, depends on $D$ and is given by Equation \eqref{eq:gamma_v}. Rain is produced when cloud exceeds a fixed threshold, $q_{precip}$, and the proportion converted, $\gamma_r$, is $1 \times 10^{-3} \text{ s}^{-1}$.

\begin{figure}[htp]
    \centering
    \includegraphics[width=\textwidth]{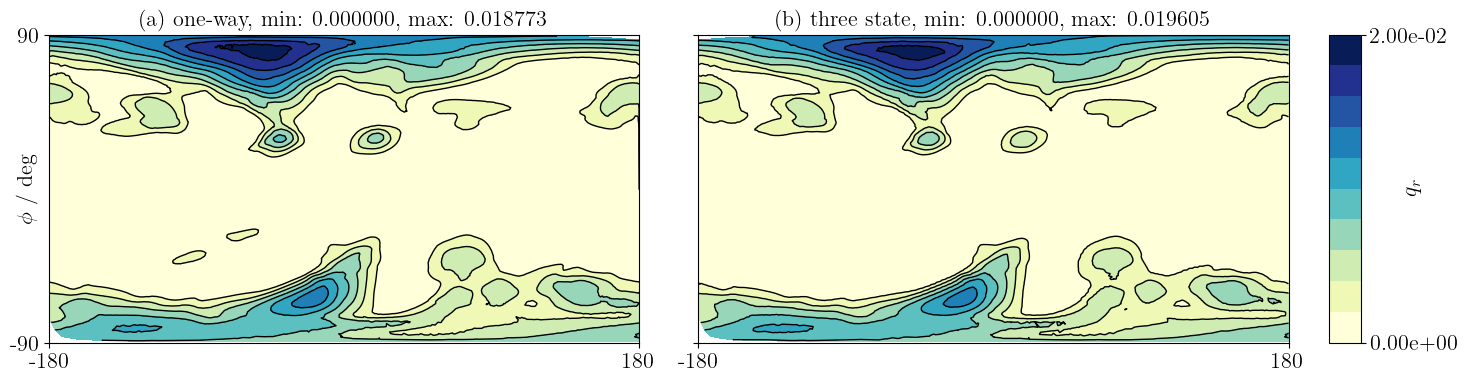}
    \caption{Comparison of the accumulated rain after 50 days in the moist convective flow over a mountain test using two different physics schemes. (a) the one-way scheme and (b) the three-state scheme.}
    \label{fig:compare_physics}
\end{figure}

Accumulated rain after 50 days using both schemes is shown in Figure \ref{fig:compare_physics}. Both schemes give a very similar rainfall pattern and rainfall amounts. The decision to use the three-state scheme in this work was motivated by our ultimate aim of models for exploring physics-dynamics coupling. The idea of a reversible conversion between vapour and cloud where cloud can evaporate in a sub-saturated atmosphere is a better representation of what happens in atmospheric models, and is an important aspect of physics-dynamics coupling. 

\section{Choice of $\beta_1$ parameter}
\begin{figure}[htp]
    \centering
    \includegraphics[width=\textwidth]{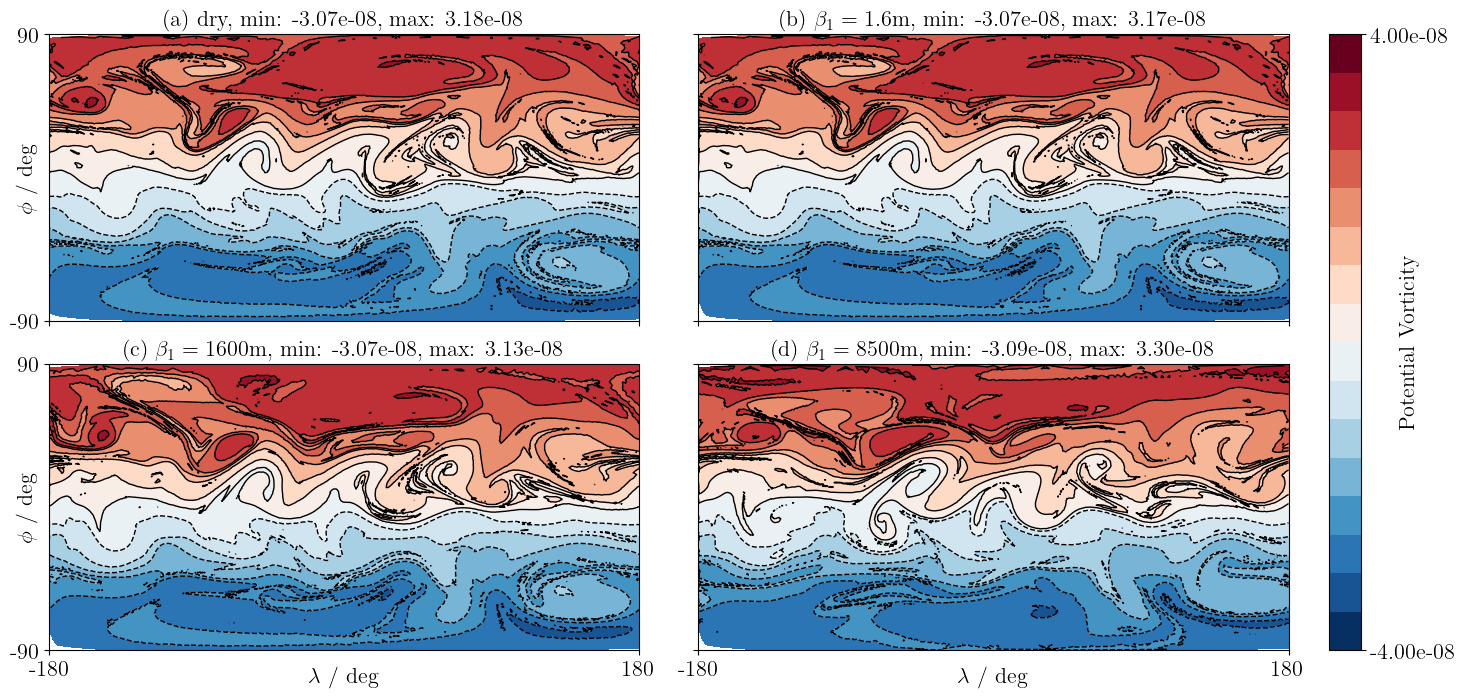}
    \caption{Comparison of the vorticity after 50 days in the standard (dry) flow over a mountain test and in the moist convective version of the same test using three different values for the $\beta_1$ parameter. (a) the dry test, (b) the moist convective test with $\beta_1 = 1.6 \text{ m}$, (c) the moist convective test with $\beta_1 = 1600 \text{ m}$, and (d) the moist convective test with $\beta_1 = 8500 \text{ m}$.}
    \label{fig:compare_beta1_mc}
\end{figure}

Our choices for $\beta_1$ and $\beta_2$ come from physically-derived values (see Appendix A in Zerroukat and Allen \cite{ZA2015} for a discussion of their pseudo-latent heat factor $L$ which is the basis for our $\beta_2$, and Equation (8) in Vallis and Penn \cite{VP2020} for an explanation for $\beta_1$). We have not done a parameter space study but we do note that varying the size of $\beta_1$ can change the impact of the convective effect. We illustrate this by comparing the potential vorticity in the standard (dry) Williamson 5 test and in the moist convective Williamson 5 test with varying $\beta_1$ values. Figure \ref{fig:compare_beta1_mc} shows the results. It is evident that the bigger the value for $\beta_1$, the bigger the difference in the large-scale dynamics as compared to the dry case. When $\beta_1$ is 1.6 m there is no discernible difference in the potential vorticity between the dry case and the moist convective case, and the difference becomes progressively bigger as $\beta_1$ is increased. As described, our choice of $\beta_1$ as 1600 m throughout the paper is motivated by a physical relationship, following Equation (8) in Vallis and Penn \cite{VP2020}. 

\begin{figure}[htp]
    \centering
    \includegraphics[width=\textwidth]{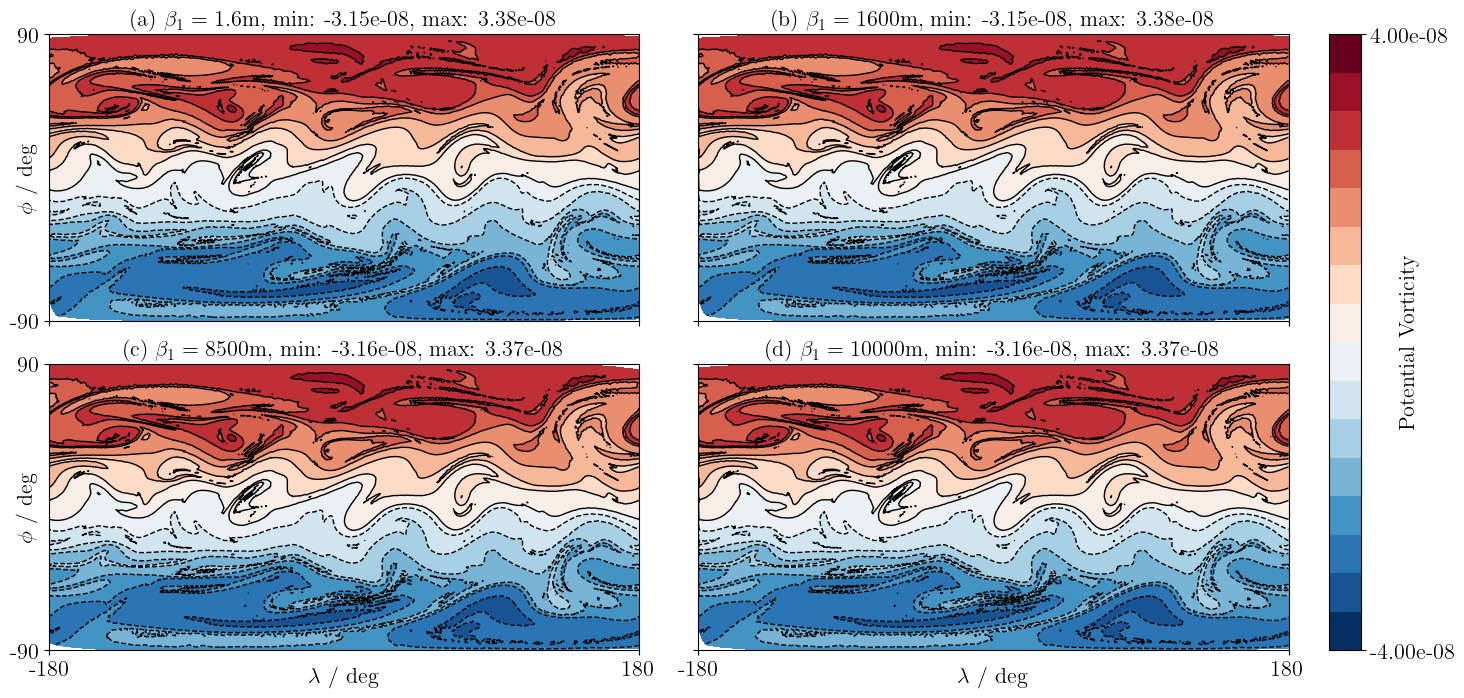}
    \caption{Comparison of the vorticity after 50 days in the moist convective thermal version of the flow over a mountain test using four different values for the $\beta_1$ parameter. (a) $\beta_1 = 1.6\text{ m}$, (b) $\beta_1 = 1600 \text{ m}$, (c) $\beta_1 = 8500 \text{ m}$, and (d) $\beta_1 = 10000 \text{ m}$. $\beta_2$ is fixed throughout at $10g \text{ m s}^{-2}$.}
    \label{fig:compare_beta1_mct}
\end{figure}

This variation in the impact of convection with varying $\beta_1$, however, is only clear in the moist convective model. In the other formulations where the buoyancy evolves in time we see very little difference in the large-scale dynamics when $\beta_1$ is varied. This is because the physics increments to depth and buoyancy are proportional to the condensation/evaporation rate, which is much bigger when the buoyancy is prescribed and does not evolve in time. In contrast, in all formulations other than the moist convective, the buoyancy evolves in a Lagrangian way, making it much less likely that a parcel of air suddenly becomes over or under-saturated. Then the physics increments to the depth field are much greater in the most convective model than in any of the others, as the condensation/evaporation is much bigger there. This is illustrated in Figure \ref{fig:compare_beta1_mct}, where we compare the potential vorticity in the moist convective thermal model with varying values of $\beta_1$. There is very little difference in the potential vorticity between the set-ups, which suggests that the physics increments from convection are small.   

\newpage
\bibliographystyle{unsrt}
\bibliography{paper}

\begin{thebibliography}{10}

\bibitem{PDCreview}
Markus Gross, Hui Wan, Philip~J Rasch, Peter~M Caldwell, David~L Williamson, Daniel Klocke, Christiane Jablonowski, Diana~R Thatcher, Nigel Wood, Mike Cullen, et~al.
\newblock Physics--dynamics coupling in weather, climate, and earth system models: Challenges and recent progress.
\newblock {\em Monthly Weather Review}, 146(11):3505--3544, 2018.

\bibitem{StechmannMajda2006}
Samuel~N Stechmann and Andrew~J Majda.
\newblock The structure of precipitation fronts for finite relaxation time.
\newblock {\em Theoretical and Computational Fluid Dynamics}, 20:377--404, 2006.

\bibitem{Gill1982}
Adrian~E. Gill.
\newblock Studies of moisture effects in simple atmospheric models: The stable case.
\newblock {\em Geophysical \& Astrophysical Fluid Dynamics}, 19(1-2):119--152, 1982.

\bibitem{Bouchut2009}
Fran{\c{c}}ois Bouchut, Julien Lambaerts, Guillaume Lapeyre, and Vladimir Zeitlin.
\newblock Fronts and nonlinear waves in a simplified shallow-water model of the atmosphere with moisture and convection.
\newblock {\em Physics of Fluids}, 21(11):116604, 2009.

\bibitem{Lambaertsetal2011baroclinic}
Julien Lambaerts, Guillaume Lapeyre, Vladimir Zeitlin, and Fran{\c{c}}ois Bouchut.
\newblock Simplified two-layer models of precipitating atmosphere and their properties.
\newblock {\em Physics of Fluids}, 23(4):046603, 2011.

\bibitem{KLZ2020b}
Alexander Kurganov, Yongle Liu, and Vladimir Zeitlin.
\newblock Moist-convective thermal rotating shallow water model.
\newblock {\em Physics of Fluids}, 32(6):066601, 2020.

\bibitem{ZA2015}
Mohamed Zerroukat and Thomas Allen.
\newblock A moist {B}oussinesq shallow water equations set for testing atmospheric models.
\newblock {\em Journal of Computational Physics}, 290:55--72, 2015.

\bibitem{Yang2021}
Da~Yang.
\newblock A shallow-water model for convective self-aggregation.
\newblock {\em Journal of the Atmospheric Sciences}, 78(2):571--582, 2021.

\bibitem{WC2014}
Michael W{\"u}rsch and George~C Craig.
\newblock A simple dynamical model of cumulus convection for data assimilation research.
\newblock {\em Meteorologische Zeitschrift}, 23(5):483--490, 2014.

\bibitem{Lambaertsetal2011jet}
Julien Lambaerts, Guillaume Lapeyre, and Vladimir Zeitlin.
\newblock Moist versus dry barotropic instability in a shallow-water model of the atmosphere with moist convection.
\newblock {\em Journal of the Atmospheric Sciences}, 68(6):1234--1252, 2011.

\bibitem{LahayeZeitlin2016}
No{\'e} Lahaye and Vladimir Zeitlin.
\newblock Understanding instabilities of tropical cyclones and their evolution with a moist convective rotating shallow-water model.
\newblock {\em Journal of the Atmospheric Sciences}, 73(2):505--523, 2016.

\bibitem{RostamiZeitlin2017}
Masoud Rostami and Vladimir Zeitlin.
\newblock Influence of condensation and latent heat release upon barotropic and baroclinic instabilities of vortices in a rotating shallow water f-plane model.
\newblock {\em Geophysical \& Astrophysical Fluid Dynamics}, 111(1):1--31, 2017.

\bibitem{RostamiZeitlin2018improved}
Masoud Rostami and Vladimir Zeitlin.
\newblock An improved moist-convective rotating shallow-water model and its application to instabilities of hurricane-like vortices.
\newblock {\em Quarterly Journal of the Royal Meteorological Society}, 144(714):1450--1462, 2018.

\bibitem{RostamiZeitlin2019}
Masoud Rostami and Vladimir Zeitlin.
\newblock Eastward-moving convection-enhanced modons in shallow water in the equatorial tangent plane.
\newblock {\em Physics of Fluids}, 31(2):021701, 2019.

\bibitem{RostamiZeitlin2020}
Masoud Rostami and Vladimir Zeitlin.
\newblock Evolution, propagation and interactions with topography of hurricane-like vortices in a moist-convective rotating shallow-water model.
\newblock {\em Journal of Fluid Mechanics}, 902:A24, 2020.

\bibitem{RostamiZeitlin2020MJO}
Masoud Rostami and Vladimir Zeitlin.
\newblock Can geostrophic adjustment of baroclinic disturbances in the tropical atmosphere explain {MJO} events?
\newblock {\em Quarterly Journal of the Royal Meteorological Society}, 146(733):3998--4013, 2020.

\bibitem{RostamiZeitlin2021}
Masoud Rostami and Vladimir Zeitlin.
\newblock Eastward-moving equatorial modons in moist-convective shallow-water models.
\newblock {\em Geophysical \& Astrophysical Fluid Dynamics}, 115(3):345--367, 2021.

\bibitem{RostamiZeitlin2022}
Masoud Rostami and Vladimir Zeitlin.
\newblock Evolution of double-eye wall hurricanes and emergence of complex tripolar end states in moist-convective rotating shallow water model.
\newblock {\em Physics of Fluids}, 34(6):066602, 2022.

\bibitem{RostamiZeitlin2022jets}
Masoud Rostami and Vladimir Zeitlin.
\newblock Instabilities of low-latitude easterly jets in the presence of moist convection and topography and related cyclogenesis, in a simple atmospheric model.
\newblock {\em Geophysical \& Astrophysical Fluid Dynamics}, 116(1):56--77, 2022.

\bibitem{VP2020}
Geoffrey~K Vallis and James Penn.
\newblock Convective organization and eastward propagating equatorial disturbances in a simple excitable system.
\newblock {\em Quarterly Journal of the Royal Meteorological Society}, 146(730):2297--2314, 2020.

\bibitem{Cotter2023}
Colin~J Cotter.
\newblock Compatible finite element methods for geophysical fluid dynamics.
\newblock {\em Acta Numerica}, 32:291--393, 2023.

\bibitem{Lawrence2018}
Bryan~N Lawrence, Michael Rezny, Reinhard Budich, Peter Bauer, J{\"o}rg Behrens, Mick Carter, Willem Deconinck, Rupert Ford, Christopher Maynard, Steven Mullerworth, et~al.
\newblock Crossing the chasm: how to develop weather and climate models for next generation computers?
\newblock {\em Geoscientific Model Development}, 11(5):1799--1821, 2018.

\bibitem{StaniforthThuburn2012}
Andrew Staniforth and John Thuburn.
\newblock Horizontal grids for global weather and climate prediction models: a review.
\newblock {\em Quarterly Journal of the Royal Meteorological Society}, 138(662):1--26, 2012.

\bibitem{GungHoCartesian}
Thomas Melvin, Tommaso Benacchio, Ben Shipway, Nigel Wood, John Thuburn, and Colin Cotter.
\newblock A mixed finite-element, finite-volume, semi-implicit discretization for atmospheric dynamics: Cartesian geometry.
\newblock {\em Quarterly Journal of the Royal Meteorological Society}, 145(724):2835--2853, 2019.

\bibitem{Williamson1992}
David~L Williamson, John~B Drake, James~J Hack, R{\"u}diger Jakob, and Paul~N Swarztrauber.
\newblock A standard test set for numerical approximations to the shallow water equations in spherical geometry.
\newblock {\em Journal of Computational Physics}, 102(1):211--224, 1992.

\bibitem{Galewsky2004}
Joseph Galewsky, Richard~K Scott, and Lorenzo~M Polvani.
\newblock An initial-value problem for testing numerical models of the global shallow-water equations.
\newblock {\em Tellus A: Dynamic Meteorology and Oceanography}, 56(5):429--440, 2004.

\bibitem{Wood2014inherently}
Nigel Wood, Andrew Staniforth, Andy White, Thomas Allen, Michail Diamantakis, Markus Gross, Thomas Melvin, Chris Smith, Simon Vosper, Mohamed Zerroukat, et~al.
\newblock An inherently mass-conserving semi-implicit semi-{L}agrangian discretization of the deep-atmosphere global non-hydrostatic equations.
\newblock {\em Quarterly Journal of the Royal Meteorological Society}, 140(682):1505--1520, 2014.

\bibitem{Zeitlinbook}
Vladimir Zeitlin.
\newblock {\em Geophysical fluid dynamics: understanding (almost) everything with rotating shallow water models}.
\newblock Oxford University Press, 2018.

\bibitem{Yang2013}
Da~Yang and Andrew~P Ingersoll.
\newblock Triggered convection, gravity waves, and the {MJO}: A shallow-water model.
\newblock {\em Journal of the atmospheric sciences}, 70(8):2476--2486, 2013.

\bibitem{Yang2014}
Da~Yang and Andrew~P Ingersoll.
\newblock A theory of the {MJO} horizontal scale.
\newblock {\em Geophysical Research Letters}, 41(3):1059--1064, 2014.

\bibitem{Kent2017}
Thomas Kent, Onno Bokhove, and Steven Tobias.
\newblock Dynamics of an idealized fluid model for investigating convective-scale data assimilation.
\newblock {\em Tellus A: Dynamic Meteorology and Oceanography}, 69(1):1369332, 2017.

\bibitem{Kent2020}
Thomas Kent, Luca Cantarello, Gordon Inverarity, Steven Tobias, and Onno Bokhove.
\newblock Idealized forecast-assimilation experiments for convective-scale numerical weather prediction.
\newblock 2020.

\bibitem{Bokhove2022}
Onno Bokhove, Luca Cantarello, and Steven Tobias.
\newblock An idealized 1$1/2$-layer isentropic model with convection and precipitation for satellite data assimilation research. part {II}: model derivation.
\newblock {\em Journal of the Atmospheric Sciences}, 79(3):875--886, 2022.

\bibitem{Vzagar2008}
Nedjeljka {\v{Z}}agar, Ad~Stoffelen, Gert-Jan Marseille, Christophe Accadia, and Peter Schl{\"u}ssel.
\newblock Impact assessment of simulated {D}oppler wind lidars with a multivariate variational assimilation in the tropics.
\newblock {\em Monthly Weather Review}, 136(7):2443--2460, 2008.

\bibitem{Vzagar2012}
Nedjeljka {\v{Z}}agar.
\newblock Multivariate data assimilation in the tropics by using equatorial waves.
\newblock {\em Pure and applied geophysics}, 169:367--379, 2012.

\bibitem{Frierson2004}
Dargan~MW Frierson, Andrew~J Majda, and Olivier~M Pauluis.
\newblock Large scale dynamics of precipitation fronts in the tropical atmosphere: A novel relaxation limit.
\newblock {\em Communications in Mathematical Sciences}, 2004.

\bibitem{Pauluis2008}
Olivier Pauluis, Dargan~MW Frierson, and Andrew~J Majda.
\newblock Precipitation fronts and the reflection and transmission of tropical disturbances.
\newblock {\em Quarterly Journal of the Royal Meteorological Society: A journal of the atmospheric sciences, applied meteorology and physical oceanography}, 134(633):913--930, 2008.

\bibitem{Lambaerts2012}
Julien Lambaerts, Guillaume Lapeyre, and Vladimir Zeitlin.
\newblock Moist versus dry baroclinic instability in a simplified two-layer atmospheric model with condensation and latent heat release.
\newblock {\em Journal of the Atmospheric Sciences}, 69(4):1405--1426, 2012.

\bibitem{Zhao2021}
Bowen Zhao, Vladimir Zeitlin, and Alexey~V Fedorov.
\newblock Equatorial modons in dry and moist-convective shallow-water systems on a rotating sphere.
\newblock {\em Journal of Fluid Mechanics}, 916:A8, 2021.

\bibitem{BettsMiller}
Alan~K Betts and Martin~J Miller.
\newblock A new convective adjustment scheme. part {II}: Single column tests using {GATE} wave, {BOMEX}, {ATEX} and {A}rctic air-mass data sets.
\newblock {\em Quarterly Journal of the Royal Meteorological Society}, 112(473):693--709, 1986.

\bibitem{KLZ2020}
Alexander Kurganov, Yongle Liu, and Vladimir Zeitlin.
\newblock Thermal versus isothermal rotating shallow water equations: comparison of dynamical processes by simulations with a novel well-balanced central-upwind scheme.
\newblock {\em Geophysical {\&} Astrophysical Fluid Dynamics}, 2020.

\bibitem{KLZ2021}
Alexander Kurganov, Yongle Liu, and Vladimir Zeitlin.
\newblock Interaction of tropical cyclone-like vortices with sea-surface temperature anomalies and topography in a simple shallow-water atmospheric model.
\newblock {\em Physics of Fluids}, 33(10):106606, 2021.

\bibitem{Rostami2022}
Masoud Rostami, Bowen Zhao, and Stefan Petri.
\newblock On the genesis and dynamics of {M}adden--{J}ulian oscillation-like structure formed by equatorial adjustment of localized heating.
\newblock {\em Quarterly Journal of the Royal Meteorological Society}, 2022.

\bibitem{SantosPeixoto2021}
Luan~F Santos and Pedro~S Peixoto.
\newblock Topography-based local spherical {V}oronoi grid refinement on classical and moist shallow-water finite-volume models.
\newblock {\em Geoscientific Model Development}, 14(11):6919--6944, 2021.

\bibitem{FergusonJablonowski2019}
Jared~O Ferguson, Christiane Jablonowski, and Hans Johansen.
\newblock Assessing adaptive mesh refinement ({AMR}) in a forced shallow-water model with moisture.
\newblock {\em Monthly Weather Review}, 147(10):3673--3692, 2019.

\bibitem{Vallis2019Book}
Geoffrey~K Vallis.
\newblock {\em Essentials of Atmospheric and Oceanic Dynamics}.
\newblock Cambridge University Press, 2019.

\bibitem{Vallis2019RainyBenard}
Geoffrey~K Vallis, Douglas~J Parker, and Steven~M Tobias.
\newblock A simple system for moist convection: the {R}ainy--{B}{\'e}nard model.
\newblock {\em Journal of Fluid Mechanics}, 862:162--199, 2019.

\bibitem{Bendall2023trilemma}
Thomas~M Bendall, Nigel Wood, John Thuburn, and Colin~J Cotter.
\newblock A solution to the trilemma of the moist {C}harney--{P}hillips staggering.
\newblock {\em Quarterly Journal of the Royal Meteorological Society}, 149(750):262--276, 2023.

\bibitem{CotterShipton2012}
Colin~J Cotter and Jemma Shipton.
\newblock Mixed finite elements for numerical weather prediction.
\newblock {\em Journal of Computational Physics}, 231(21):7076--7091, 2012.

\bibitem{ShiptonCotterGibson}
Jemma Shipton, Thomas~H Gibson, and Colin~J Cotter.
\newblock Higher-order compatible finite element schemes for the nonlinear rotating shallow water equations on the sphere.
\newblock {\em Journal of Computational Physics}, 375:1121--1137, 2018.

\bibitem{Natale2016}
Andrea Natale, Jemma Shipton, and Colin~J Cotter.
\newblock Compatible finite element spaces for geophysical fluid dynamics.
\newblock {\em Dynamics and Statistics of the Climate System}, 1(1):dzw005, 2016.

\bibitem{Bendall2020compatible}
Thomas~M Bendall, Thomas~H Gibson, Jemma Shipton, Colin~J Cotter, and Ben Shipway.
\newblock A compatible finite-element discretisation for the moist compressible {E}uler equations.
\newblock {\em Quarterly Journal of the Royal Meteorological Society}, 146(732):3187--3205, 2020.

\bibitem{CotterShipton2023}
Colin~J Cotter and Jemma Shipton.
\newblock A compatible finite element discretisation for the nonhydrostatic vertical slice equations.
\newblock {\em GEM-International Journal on Geomathematics}, 14(1):25, 2023.

\bibitem{NataleCotter2018}
Andrea Natale and Colin~J Cotter.
\newblock A variational finite-element discretization approach for perfect incompressible fluids.
\newblock {\em IMA Journal of Numerical Analysis}, 38(3):1388--1419, 2018.

\bibitem{Gibson2020slate}
Thomas~H Gibson, Lawrence Mitchell, David~A Ham, and Colin~J Cotter.
\newblock Slate: extending {F}iredrake's domain-specific abstraction to hybridized solvers for geoscience and beyond.
\newblock {\em Geoscientific Model Development}, 13(2):735--761, 2020.

\bibitem{Bendall2019recovered}
Thomas~M Bendall, Colin~J Cotter, and Jemma Shipton.
\newblock The ‘recovered space’ advection scheme for lowest-order compatible finite element methods.
\newblock {\em Journal of Computational Physics}, 390:342--358, 2019.

\bibitem{Kuzmin2010vertex}
Dmitri Kuzmin.
\newblock A vertex-based hierarchical slope limiter for p-adaptive discontinuous {G}alerkin methods.
\newblock {\em Journal of Computational and Applied Mathematics}, 233(12):3077--3085, 2010.

\bibitem{CotterKuzminlimiters}
Colin~J Cotter and Dmitri Kuzmin.
\newblock Embedded discontinuous {G}alerkin transport schemes with localised limiters.
\newblock {\em Journal of Computational Physics}, 311:363--373, 2016.

\bibitem{Firedrakebook}
Thomas~H Gibson, Andrew~TT McRae, Colin~J Cotter, Lawrence Mitchell, and David~A Ham.
\newblock {\em Compatible Finite Element Methods for Geophysical Flows: Automation and Implementation Using Firedrake}.
\newblock Springer Nature, 2019.

\bibitem{FiredrakeUserManual}
David~A. Ham, Paul H.~J. Kelly, Lawrence Mitchell, Colin~J. Cotter, Robert~C. Kirby, Koki Sagiyama, Nacime Bouziani, Sophia Vorderwuelbecke, Thomas~J. Gregory, Jack Betteridge, Daniel~R. Shapero, Reuben~W. Nixon-Hill, Connor~J. Ward, Patrick~E. Farrell, Pablo~D. Brubeck, India Marsden, Thomas~H. Gibson, Miklós Homolya, Tianjiao Sun, Andrew T.~T. McRae, Fabio Luporini, Alastair Gregory, Michael Lange, Simon~W. Funke, Florian Rathgeber, Gheorghe-Teodor Bercea, and Graham~R. Markall.
\newblock {\em Firedrake User Manual}.
\newblock Imperial College London and University of Oxford and Baylor University and University of Washington, first edition, 5 2023.

\bibitem{UFL}
Martin~S Aln{\ae}s, Anders Logg, Kristian~B {\O}lgaard, Marie~E Rognes, and Garth~N Wells.
\newblock Unified form language: A domain-specific language for weak formulations of partial differential equations.
\newblock {\em ACM Transactions on Mathematical Software (TOMS)}, 40(2):1--37, 2014.

\bibitem{PETSc}
Satish Balay, Shrirang Abhyankar, Mark Adams, Jed Brown, Peter Brune, Kris Buschelman, Lisandro Dalcin, Alp Dener, Victor Eijkhout, William Gropp, et~al.
\newblock Petsc users manual.
\newblock 2019.

\end{thebibliography}
\end{document}